\documentclass[11pt]{article}

\usepackage{epsfig,graphics}
\usepackage{amssymb,amsmath,amsthm,bm}

\newtheorem{theo}{Theorem}
\newtheorem{lem}{Lemma}
\newtheorem{pro}{Proposition}

\newtheorem{rem}{Remark}
\newtheorem{mdef}{Definition}

\newtheorem{alg}{Algorithm}

\def\RR{{\mathbb R}}

\def\R{{\mathbb R}}

\def\gga{\gamma}

\def\gb{\beta}

\def\gs{\sigma}
\def\gl{\lambda}
\def\wt{\widetilde}
\def\n{\noindent}

\def\gUp{\Upsilon}

\def\ep{\epsilon}
\def\vep{\varepsilon}

\def\gt{\triangle}

\def\gp{{\prime}}
\def\b0{{\bf 0}}
\def\1{{\bf 1}}

\def\cE{\mathcal {E}}
\def\cG{\mathcal G}

\def\cP{\mathcal P}
\def\cQ{\mathcal Q}

\def\Err{{\rm Res}}

\def\tErr{{\rm Err}}

\def\Bd{{\rm Bd}}

\def\wh{\widehat}
\def\wt{\widetilde}

\usepackage{color}

\makeatletter
\def\widebreve{\mathpalette\wide@breve}
\def\wide@breve#1#2{\sbox\z@{$#1#2$}%
     \mathop{\vbox{\m@th\ialign{##\crcr
\kern0.08em\brevefill#1{0.8\wd\z@}\crcr\noalign{\nointerlineskip}%
                    $\hss#1#2\hss$\crcr}}}\limits}
\def\brevefill#1#2{$\m@th\sbox\tw@{$#1($}%
  \hss\resizebox{#2}{\wd\tw@}{\rotatebox[origin=c]{90}{\upshape(}}\hss$}
\makeatletter
\def\wb{\widebreve}

\begin{document}

\title{Analysis of  a Direct Separation Method
Based on Adaptive Chirplet Transform for Signals  with Crossover Instantaneous Frequencies\thanks{This work is partially supported by 
the ARO under Grant $\sharp$ W911NF2110218,
the Simons Foundation under grant $\sharp$ 353185, and the National Natural Science Foundation of China 
under Grants $\sharp$ 62071349, $\sharp$ 61972265, $\sharp$ 11871348 and $\sharp$ U21A20455,.
}
\\
{\small (submitted 
on Jan 30, 2022)}}
\author{Charles K. Chui${}^{1}$, 
 Qingtang Jiang${}^{2}$, Lin Li${}^{3}$, and Jian Lu${}^{4}$
  }
\date{}

\maketitle

\vskip -0.5cm


\bigskip
{\small 1. Menlo Park residence, CA 94025, USA. e-mail: ckchui@stanford.edu.
}

{\small 2. Department of Mathematics \& Statistics, 
University of Missouri-St. Louis,} 

\quad {\small St. Louis,  MO 63121, USA.} 

{\small 3. School of Electronic Engineering, Xidian University, Xi${}'$an 710071, China.}  

{\small 4. Shenzhen Key Laboratory of Advanced Machine Learning and Applications,}  

\quad {\small College of Mathematics \& Statistics, Shenzhen University, Shenzhen 518060, China.} 

\begin{abstract}
In many applications, it is necessary to retrieve the sub-signal building blocks of a multi-component signal, which is usually non-stationary in real-world and real-life applications. Empirical mode decomposition (EMD), synchrosqueezing transform (SST), signal separation operation (SSO), and iterative filtering decomposition (IFD) have been proposed and developed for this purpose. However, these computational methods are restricted by the specification of well-separation of the sub-signal frequency curves for multi-component signals. On the other hand, the chirplet transform-based signal separation scheme (CT3S) that extends SSO from the two-dimensional ``time-frequency" plane to the three-dimensional ``time-frequency-chirp rate" space was recently proposed in our recent work to remove the frequency-separation specification, and thereby allowing ``frequency crossing". The main objective of this present paper is to carry out an in-depth error analysis study of instantaneous frequency estimation and component recovery for the CT3S method.
\end{abstract}



\maketitle

\section{Introduction}

Many real-world signals consist of sub-signal building blocks (also called components or modes of the signal). To analyze such multi-component signals, it is necessary to recover the signal components.  The empirical mode decomposition (EMD) scheme \cite{Huang98}, the synchrosqueezing transform (SST) \cite{Daub_Lu_Wu11,Thakur_Wu11}, the signal separation operator (SSO) 
 method \cite{Chui_M16, Chui_Mhaskar_Walt}, and the iterative filtering decomposition (IFD) scheme \cite{HM_Zhou09, HM_Zhou16} are among the existing computational methods developed for this purpose. In particular, the SSO scheme is a direct time-frequency method without the need of computing a ``reference frequency curve" for the SST approach and repeated iterations for the EMD and IFD computational schemes. For this reason, based on the SSO method, we proposed the ``chirplet transform-based signal separation scheme" in our recent paper \cite{LHJC22}. 
 
 Let $x(t)$ be a nonstationary multi-component signal defined by 
\begin{equation}
\label{MAHM}
x(t)=A_{0}(t)+\sum_{k=1}^K x_k(t)=A_0(t)+\sum_{k=1}^K A_k(t) e^{i2\pi \phi_k(t)},   
\end{equation}
where $A_k(t), \phi_k'(t)$ are positive functions, such that for each $k$, $0\le k\le K$, as a function of $t$, $A_k(t)$ ``changes very slowly". 
Let $V_x(t, \eta)$ be the  (modified)  short-time Fourier transform (STFT) of $x(t)$ defined by 
\begin{equation}
\label{def_STFT}
V_x(t, \eta):=\int_{\RR} x(\tau) g(\tau-t) e^{-i2\pi \eta(\tau-t)}d\tau,
\end{equation}
where $g(t)$ is a window function. Under certain conditions, the ridges $\wh \eta_k(t)$ of the spectrogram $|V_x(t, \eta,\lambda)|$ provide an estimate of the instantaneous frequencies (IFs) $\phi^\gp_k(t)$, so that the SSO method can be applied to approximately reconstruct $x_k(t)$ by 
$$
x_k(t)\approx V_x(t, \wh \eta_k(t)). 
$$
Observe that in comparison with the SST approach, SSO does not require the ``squeezing" operation to compute some ``reference frequency" as required by SST, but directly retrieves the components $x_k(t)$ simply by replacing $\eta$ in $V_x(t, \eta)$ by the time-frequency ridges $\wh \eta_k(t)$. Hence, we refer this ``modified SSO" method as a direct time-frequency approach.  An improvement of the SSO method based on linear chirp local approximation was  recently introduced \cite{LChuiJ20, 
CJLL20_adpSTFT}. In this regard, we also mention the continuous wavelet transform (CWT)-based SSO signal separation methods introduced and developed in \cite{Chui_Han20, CJLL21_adpCWT}. 

As point out above that the application of the EMD, SST, SSO or IFD methods is restricted to multi-component signals with well-separated frequency curves, namely: $\phi^\gp_k(t)$ are are required to satisfy 
\begin{equation}
 \label{def_sep_cond}
 \big| \phi'_k(t)-\phi'_{\ell}(t)\big| \ge 2\gt, \; 2\le k, \ell\le K, k\not=\ell;
\end{equation}
or \begin{equation}
\label{freq_resolution}
\frac{\big|\phi_k'(t)-\phi'_\ell(t)\big|}{\phi_k'(t)+\phi'_\ell(t)}\ge \gt, 
\; 2\le k, \ell\le K,   k\not=\ell 
\end{equation}
for some $\gt>0$ and all  $t\in \R$. However, in applications there are multi-component signals, such as radar echoes with micro-Doppler effect (micro motion), that overlap in the time-frequency plane  \cite{radar_basic_2016, Stankovic_compressive_sense_2013}. Here, we say $x_{k-1}(t)$ and $x_k(t)$ in \eqref{MAHM} overlap in time-frequency plane or they have crossover frequencies at $t=t_0$, 
if $\phi_{k-1}^\gp(t_{0})=\phi_{k}^\gp(t_{0})$. While none of the aforementioned methods can be applied to retrieve the modes of a multi-component signal overlapping in the time-frequency plane accurately, the chirplet transform-based signal separation scheme (CT3S for short) was proposed in \cite{LHJC22} to meet this challenge. To summarise the proposed method in \cite{LHJC22}, let $\cQ_x(t, \eta, \gl)$ be the (adaptive) chirplet transform (CT) of $x(t)$ defined by 
\begin{eqnarray}
 \label{def_ACT}
 \cQ_x(t,\eta,\lambda)\hskip -0.6cm && := \int_{\RR} x(\tau) \frac 1 {\gs(t)} g\big(\frac {\tau-t} {\gs(t)}\big) e^{-i2\pi\eta(\tau-t) -i\pi \lambda (\tau-t)^2}
 d\tau\\
\nonumber && = \int_{\RR} x(t+\tau) 
\frac 1 {\gs(t)} g\big(\frac \tau {\gs(t)}\big) e^{-i2\pi\eta\tau -i\pi \lambda \tau^2}
d\tau,
 \end{eqnarray}
where $g(t)$ is a window function and $\gs(t)$ is a positive function. In the published literature, $\cQ_x$ is also called the localized polynomial Fourier transform (of order 2),  see, for example, \cite{Bi_Stankovic11}. The CT3S is described as follows. Let $x(t)$, as given in  \eqref{MAHM} with $\phi_0(t): \equiv 0$, satisfy   
\begin{equation}\label{def_sep_cond_cros}
	|\phi'_{k}(t)-\phi'_{\ell}(t)|+\rho \; |\phi''_{k}(t)-\phi''_{\ell}(t)| \ge 2 \gt, \; 0\le \ell, k\le K, \ell\not=k, 
	\end{equation}
where $\gt>0, \rho>0$. 
For a fixed $t$, and a positive number  $\wt \ep_1$, 
consider the sets: 
\begin{equation}
\label{def_cGk}
\begin{array}{l}
\cG_t:=\{(\eta, \gl): \; |\cQ_x(t, \eta, \gl)|>\wt \ep_1\}; \\
\cG_{t, k}:=\{(\eta, \gl) \in \cG_t: \;  
|\eta-\phi_k^\gp(t)|+ \rho \; |\gl-\phi''_k(t)|< \gt\}.  
\end{array}
\end{equation}
Under certain conditions (see Theorem 1 in the next section),  $\cG_t$ is a disjoint union of non-empty sets  $\cG_{t, k}$, where $0\le k\le K$.  In the following algorithm, we need the notations:
\begin{equation}
\label{def_max_eta}
\wh \eta_0:=0, \; \wh \gl_0:=0, \; (\wh \eta_k, \wh \gl_k) =(\wh \eta_k(t), \wh \gl_k(t)):={\rm argmax}_{(\eta, \gl) \in\mathcal{G}_{t, k}  }|\cQ_x(t,\eta, \gl)|, \; k=1, \cdots, K.
\end{equation}

 \begin{alg} {\bf (Chirplet transform-based signal separation scheme (CT3S))} \; For  $x(t)$ defined in \eqref{MAHM}, do the following 
\begin{itemize}
\item[] {\bf Step 1.} Calculate  $\wh \eta_k(t)$ and $\wh \gl_k(t)$ by \eqref{def_max_eta}.
\item[] {\bf Step 2.} Obtain IF and chirp rate estimates by 
\begin{equation}
\label{IF_estimate}
\phi^\gp_k(t) \approx \wh \eta_k(t), \quad 
\phi^{\gp\gp}_k(t) \approx \wh \gl_k(t), \; 1\le k \le K. 
\end{equation} 
\item[] {\bf Step 3.} Obtain the recovered $k$-th component by 
\begin{equation}
\label{comp_recover}
x_k(t)\approx \cQ_x\big(t, \wh \eta_k(t), \wh \gl_k(t)\big), \; 0\le k \le K.  
\end{equation}
\hfill $\blacksquare$ 
\end{itemize}
\end{alg}

\begin{figure}[th]
	\centering
	\begin{tabular}{cc}
	\resizebox {2.2in}{1.5in} {\includegraphics{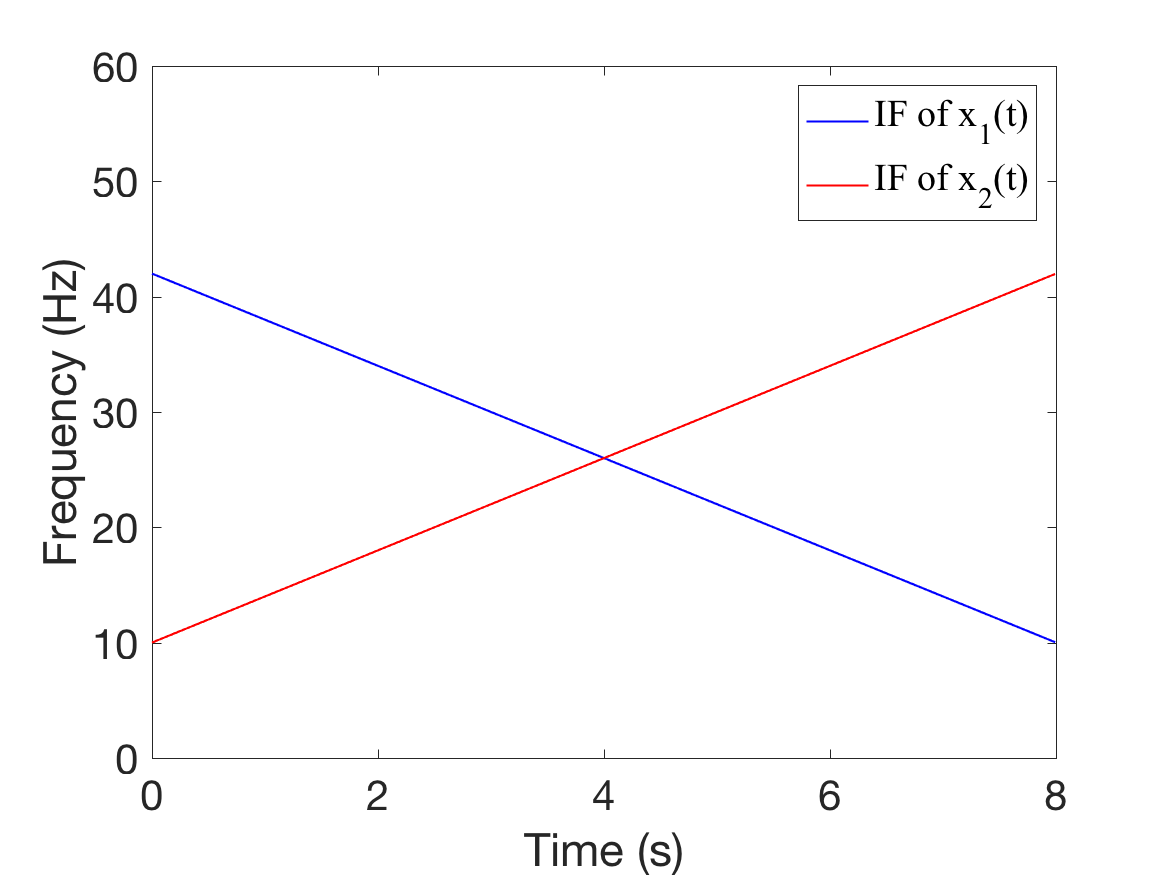}}
	&
\resizebox {2.2in}{1.5in}{\includegraphics{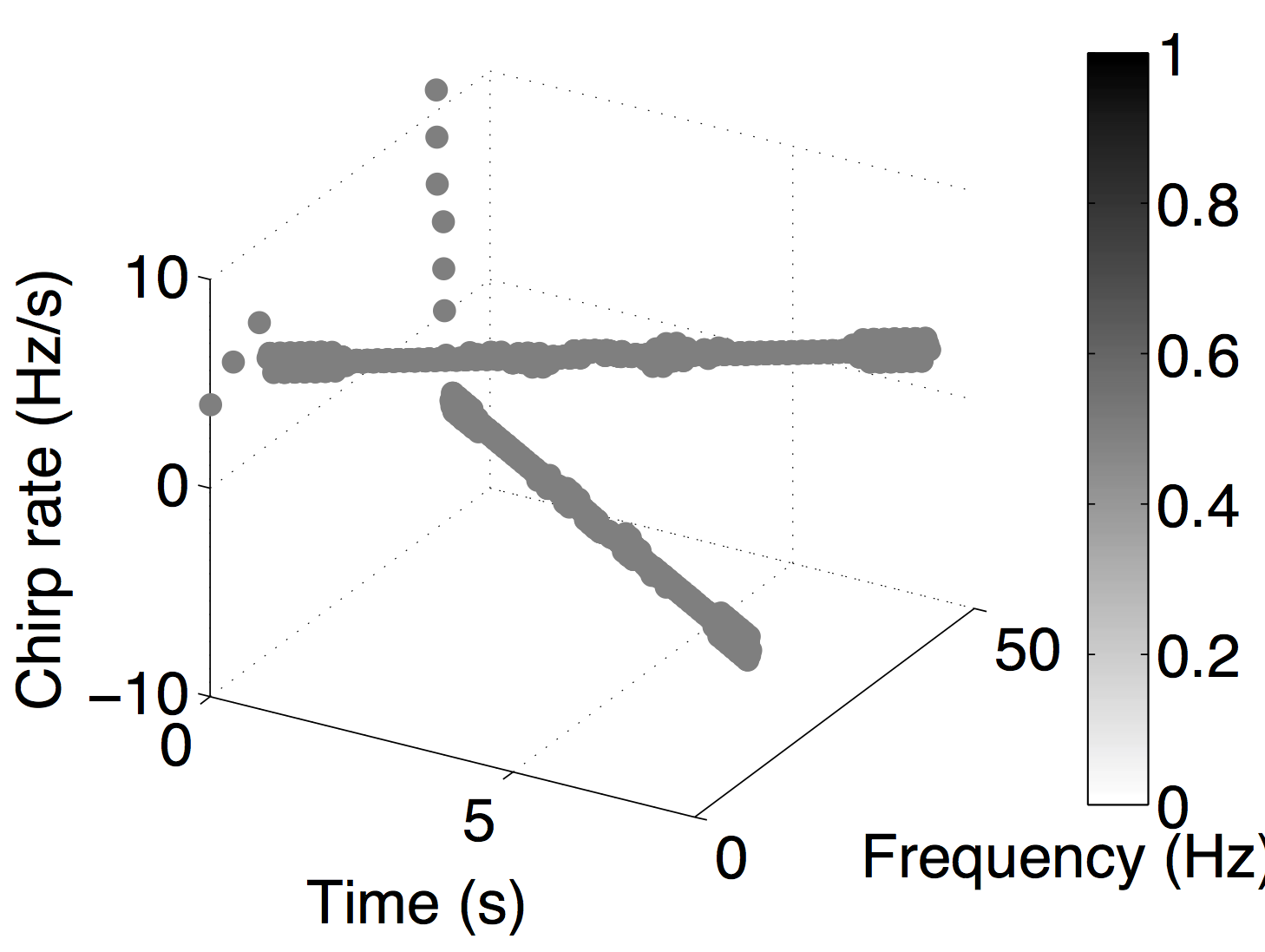}}
	\end{tabular}
	\vskip -0.3cm 
	\caption{\small  Left: IFs of $x_1, x_2$; Right:  
	Ridges $(t, \wh \eta_1(t), \wh \gl_1(t))$  and $(t, \wh \eta_2(t), \wh \gl_2(t))$. 
	} 
	\label{figure:3DTFCr12}
\end{figure}

It is worth noting that more general notions of the chirplet tansform (CT) were already introduced in \cite{Mann95}, and such CT was used for IF estimation and mode retrieval of multi-component signals in the recent literature, such as ``multi-synchrosqueezing CT" for IF estimation in \cite{Zhu et al 2019}; CT-based joint estimations of IFs and chirp rates in \cite{Yang_CT2020}, and synchrosqueezed CT for IF estimation and mode retrieval in \cite{Wu_SCT2021}. The interested reader is referred to \cite{LHJC22, Wu_SCT2021} for other references on CT-based IF estimation and/or component recovery. However, the method proposed in our paper \cite{LHJC22} differs from the others, in that
``ridges" and the formula \eqref{comp_recover} were introduced to recover signal components of multi-component signals. Next we provide an example to illustrate why the idea of our proposed direct method of CT-based 3-D transform works in recovering multi-component signal modes, even if the modes overlap in the time-frequency plane. Let  
\begin{equation}
\label{x_example1}
x(t)=x_1(t)+x_2(t), \; x_1(t)=e^{i2\pi (42t-2t^2)}, \; 
 x_2(t)=e^{i2\pi (10t+2t^2)}, \; 0\le t\le 8.
 \end{equation} 
Then the IFs of $x_1(t)$ and $x_2(t)$ are $\phi^\gp_1(t)=42-4t$,  $\phi^\gp_2(t)=10+4t$, respectively, with crossover at $t=4$. 
Observe that although $x_1(t)$  and $x_2(t)$ overlap in the time-frequency plane,  the local maximum points $(\wh \eta_1(t), \wh \gl_1(t))$ and $(\wh \eta_2(t), \wh \gl_2(t))$ of $|\cQ_{x}(t, \eta, \gl)|$ (in variables $\eta$ and $\gl$), that correspond to $x_1$ and $x_2$, are well separated in the three dimensional space. See the right panel of Figure \ref{figure:3DTFCr12} for the ridges $(t, \wh \eta_1(t), \wh \gl_1(t))$ 
and $(t, \wh \eta_2(t), \wh \gl_2(t))$, which essentially lie in two different planes $(t, \eta, -4)$ and $(t, \eta, 4)$, respectively.
Indeed, both $x_1$ and $x_2$ can be recovered by  
$$
x_1(t)\approx \cQ_x(t, \wh \eta_1(t), \wh \gl_1(t)), \quad x_2(t)\approx \cQ_x(t, \wh \eta_2(t), \wh \gl_2(t)).
$$     	   

In \cite{LHJC22}, error bounds for IF estimation and component recovery in \eqref{IF_estimate} and \eqref{comp_recover}, respectively, are provided. 
It is important to point out that in order to reduce the errors, it is necessary to increase $\gs(t)$. Unfortunately, this results in poor time resolution and causes other problems for non-stationary signals with finite time duration. In our recent paper \cite{CJLL21_adpCWT_cross} the time-scale-chirprate transform, based on adaptive continuous wavelet transform (CWT), is introduced  to recover modes and our analysis there gets around such problems. In this present paper, the idea in \cite{CJLL21_adpCWT_cross} for the CWT is adopted to the STFT to carry out an in-depth study on error estimates of IF estimation and on component recovery by CT3S. The error bounds so obtained depend on certain properties of the window function $g$.  

This paper is organized as follows. In Section 2, we will establish error bounds for IF estimation and component recovery when CT3S is applied. More explicit and detailed formulation of the error bounds will be derived for the Gaussian window function $g$ in Section 3. Two numerical experiments for illustration of the general theorem are to be presented in the final section. 

\section {Signal separation by adaptive chirprate transform}

For convenience of our presentation, we reformulate the trend $A_0(t)$ in \eqref{MAHM} as $x_0(t)=A_0(t) e^{i2\pi \phi_0(t)}$ with $\phi_0(t)=0$. 
As in \cite{CJLL21_adpCWT_cross}, 
for $\ep_1>0$ and $\ep_3>0$, we let $\cE_{\ep_1,\ep_3}$ denote the set of such signals $x(t)$ in the form of \eqref{MAHM} that $A_k(t), \phi_k(t)$ satisfy that $A_k(t) \in L_\infty(\RR), \; A_k(t)>0, \phi_k(t)\in C^3(\RR), \inf_{t\in \RR} \phi_k'(t)>0, \sup_{t\in \RR} \phi_k'(t)<\infty$,  and 
\begin{eqnarray}
\label{cond_A}
&&	|A_k(t+\tau)-A_k(t)|\leq \vep_1 |\tau|A_{k}(t),~~t\in \RR, \; k=0,\cdots, K, \\
\label{cond_phi}
&&	|\phi'''_{k}(t)| \leq \vep_3,  ~~t\in \RR, \; k=1, \cdots, K. 
	\end{eqnarray}

In this paper a window function $g$ is a continuous function in $L_2(\R)\cap L_1(\R)$ with $\int_\RR g(t) dt =1$ and certain rate of decay at $\infty$ to be specified. 
For such window functions $g$, we consider 
\begin{equation}
\label{def_PFT}
\wb g (\eta, \gl):=\int_{\RR} g(\tau) e^{-i2\pi\eta\tau-i\pi \lambda \tau^2}d\tau, 
\end{equation}
to be called the ``polynomial Fourier transform (of order 2)" corresponding to $g$, as coined in  \cite{Bi_Stankovic11}. 
Also, recall a function $s(t)$ is said to be a linear frequency modulation (LFM)  signal (or a linear chirp), if
\begin{equation*}
s(t)=A e^{i2\pi \phi(t)}:=A e^{i2\pi (ct +\frac 12 r t^2)},   
\end{equation*}
for some $c, r\in \RR$, with $r\not=0$, and $A>0$. As in \cite{LChuiJ20, CJLL20_adpSTFT}, 
we will apply LFMs to approximate each sub-signal $x_k(t)$ locally at each time instant $t$ of the multi-component signal $x(t)$, in that  
$$
x_k(t+\tau)=x_k(t)e^{i2\pi (\phi^\gp_k(t)\tau+\frac 12\phi^{\gp \gp}_k(t)\tau^2) }+ x_{{\rm r}, k}(t, \tau),
$$
for small $\tau$, where 
\begin{eqnarray}
\nonumber x_{{\rm r}, k}(t, \tau)\hskip -0.6cm && :=x_k(t+\tau)-x_k(t)e^{i2\pi (\phi^\gp_k(t)\tau+\frac 12\phi^{\gp \gp}_k(t)\tau^2) }\\
\label{xrk}&&= (A_k(t+\tau)-A_k(t))e^{i2\pi \phi_k(t+\tau)}
\\
\nonumber &&\hskip 1cm 
+x_k(t)e^{i2\pi (\phi^\gp_k(t)\tau+\frac 12\phi^{\gp \gp}_k(t)\tau^2) }
\big(
e^{i2\pi (\phi_k(t+\tau)-\phi_k(t)-\phi_k^\gp(t) \tau- \frac 12\phi^{\gp \gp}_k(t)\tau^2)}-1\big).
\end{eqnarray}
Thus, we have 
$$
x(t+\tau)=x_{\rm m}(t, \tau)+ x_{\rm r}(t, \tau),
$$
where 
$$
x_{\rm m}(t, \tau):=\sum_{k=0}^K x_k(t)e^{i2\pi (\phi^\gp_k(t)\tau+\frac 12\phi^{\gp \gp}_k(t)\tau^2) }, \; 
x_{\rm r}(t, \tau):=\sum_{k=0}^K x_{{\rm r}, k}(t, \tau). 
$$
By setting  
\begin{equation}
	\label{def_MSSO_m1}
	 \cP_x(t,\eta,\lambda) :=  
 \int_{\RR} x_{\rm m}(t, \tau)  \frac 1 {\gs(t)} g\big(\frac \tau {\gs(t)}\big)  e^{-i2\pi\tau\eta-i\pi \lambda \tau^2} d\tau, 
	\end{equation}
we have 
\begin{equation}
\label{def_MSSO_m2}
 \cP_x(t,\eta,\lambda) =\sum_{k=0}^K x_k(t) \wb g\big( \gs(\eta -\phi^\gp_k(t)), \gs^2(\gl - \phi^{\gp \gp}_k(t))\big). 
	\end{equation}
Note that $\gs$ in \eqref{def_MSSO_m2} is actually $\gs(t)$. Here and there-after, we write
$\gs = \gs(t)$ for simplicity. In the following, we denote 
 \begin{equation}\label{def_M_u}
\mu=\mu(t):= \min_{0\leq k\leq K} A_{k}(t), ~~ M=M(t):=\sum_{k=0}^K A_{k}(t).  
\end{equation}

 In the following, we derive an error bound for $|\cQ_x(t, \eta, \gl)-\cP_x(t, \eta, \gl)|$. 
\begin{lem}\label{lem1} 
$\cQ_x(t, \eta, \gl)$ and $\cP_x(t, \eta, \gl) $of an $x(t)$ in $\cE_{\ep_1,\ep_3}$ for some $\ep_1>0$ and $\ep_3>0$ satisfy  
\begin{equation}
\label{S_R_error}
\big|\cQ_x(t, \eta, \gl)-\cP_x(t, \eta, \gl)\big|\le \Pi(t)M(t), 
\end{equation}   
where 
\begin{eqnarray}
\label{def_Pi0}
&&\Pi(t):=\vep_1 I_1 \gs(t) + \frac {\pi}3\vep_3  I_3 \gs^3(t), \\
\label{def_In}
 && I_n:=\int_\RR |g(t) t^n| dt, \; n=1, 2, \cdots. 
\end{eqnarray}
\end{lem}

\n {\it Proof}. \; By applying the assumptions \eqref{cond_A} and \eqref{cond_phi}, it follows from \eqref{xrk} that 
\begin{eqnarray*}
&&|x(t+\tau)-x_{\rm m}(t, \tau)|=|x_{\rm r}(t, \tau)|\\
&& \le \sum_{k=0}^K\Big\{ |A_k(t+\tau)-A_k(t)| + A_k(t) ~\big| i2\pi \big(\phi_k(t+\tau)-\phi_k(t)-\phi_k^\gp(t) \tau- \frac 12\phi^{\gp \gp}_k(t)\tau^2\big) \big|\Big\}\\
&& \le \sum_{k=0}^K\Big\{ A_k(t) \vep_1|\tau|  + A_k(t)  2\pi \sup_{\xi \in \RR}  \frac 16 \big |\phi^{\gp\gp\gp}_k(\xi) \tau^3 \big|\Big\}\\
&&\le  M(t) \vep_1 |\tau|  + M(t) \frac \pi 3 \vep_3  |\tau|^3. 
\end{eqnarray*}
Therefore, we have
\begin{eqnarray*}
&&\big|\cQ_x(t, \eta, \gl)-\cP_x(t, \eta, \gl)\big| =
\Big| \int_\RR (x(t+\tau)-x_{\rm m}(t, \tau))\frac 1{\gs(t)}g(\frac \tau{\gs(t)}) e^{-i2\pi (\eta \tau+\frac 12 \gl \tau^2)}d\tau\Big|\\
&&\qquad  \le\int_\RR  M(t) \big(\vep_1 |\tau|  + \frac \pi 3 \vep_3  |\tau|^3\big)
|\frac 1{\gs(t)}g(\frac \tau{\gs(t)}) |d\tau\\
&&\qquad =M(t)\big(\vep_1 I_1 \gs(t) + \frac {\pi}3\vep_3  I_3 \gs^3(t)\big), 
\end{eqnarray*}
so that \eqref{S_R_error} holds as desired. 
\hfill $\blacksquare$ 

Let $\rho> 0$ be the multiplication constant in \eqref{def_sep_cond_cros}. Denote 
\begin{equation}
\label{def_Zk}
Z_k:=\{(t, \eta, \gl): 
	|\eta-\phi'_k(t)|+\rho \; |\gl-\phi''_k(t)| <  \gt, \; t\in \RR\}.   
\end{equation}
Let $\Upsilon(t), \Upsilon_{\ell, k}(t)$ with $\Upsilon(t)\ge \Upsilon_{\ell, k}(t)$ for $k\not=\ell$  
be some functions that satisfy 
 \begin{equation}
 \label{def_upper_bounds}
 \begin{array}{l}
 \sup_{(\eta, \gl)\not \in \cup_{k=0}^K Z_k}\big|\wb g\big( \gs(\eta -\phi^\gp_k(t)), \gs^2(\gl - \phi^{\gp \gp}_k(t))\big)\big|\le \gUp(t), \\
\sup_{(\eta, \gl)\in Z_\ell}\big|\wb g\big( \gs(\eta -\phi^\gp_k(t)), \gs^2(\gl - \phi^{\gp \gp}_k(t))\big)\big | \le \gUp_{\ell,  k}(t). 
 \end{array}
 \end{equation}
For the Gaussian window function  $g$ defined by 
 \begin{equation}
\label{def_g}
g(t)=\frac 1{\sqrt {2\pi}} \; e^{-\frac {t^2}2},   
\end{equation}
in \S3, we will derive the explicit expressions for the quantities $\Upsilon(t)$ and $\Upsilon_{\ell, k}(t)$. 
 
\bigskip 
To establish our main theorem,  we also need the following lemma, where $\sum_{k\not= \ell}$ denotes $\sum_{\{k: ~ k\not= \ell, 0\le k\le K\}}$.
\begin{lem}
Let $x(t)\in \cE_{\ep_1,\ep_3}$ for some $\ep_1, \ep_3>0$.
Then for any  $(\eta, \gl)\in \cG_{t, \ell}$,
\begin{eqnarray}
\label{U_with_ell}
&&\big|\cQ_x(t, \eta, \gl)-x_\ell(t) \wb g\big( \gs(\eta -\phi^\gp_\ell(t)), \gs^2(\gl - \phi^{\gp \gp}_\ell(t))\big) \big|\le \Err_\ell(t),    
\end{eqnarray}   
where 
\begin{equation}
\label{def_Err}
\Err_\ell(t):=M(t)\Pi(t)+\sum_{k\not= \ell} A_k(t) \gUp_{\ell, k}(t)
\end{equation}
and $\Pi(t)$ is given in \eqref{def_Pi0}. 
\end{lem}

\n {\it Proof}. \; Let $\cP_x(t, \eta, \gl)$ be the quantity defined by \eqref{def_MSSO_m2}. Then 
for any $(\eta, \gl)\in \cG_{t, \ell}$, 
\begin{eqnarray}
\nonumber &&	\big|\cP_x(t, \eta, \gl)-x_\ell(t) \wb g\big( \gs(\eta -\phi^\gp_\ell(t)), \gs^2(\gl - \phi^{\gp \gp}_\ell(t))\big) \big| \\
\nonumber && =\Big|\sum_{k\not= \ell } x_k(t) \wb g\big( \gs(\eta -\phi^\gp_k(t)), \gs^2(\gl - \phi^{\gp \gp}_k(t))\big)\Big|\\
\label{middel_step} && \le \sum_{k\not= \ell } A_k(t)  \Big| \wb g\big( \gs(\eta -\phi^\gp_k(t)), \gs^2(\gl - \phi^{\gp \gp}_k(t))\big)\Big| 
\le \sum_{k\not= \ell } A_k(t)  \gUp_{\ell, k}(t).  
\end{eqnarray}
Hence 
\begin{eqnarray*}
&&	\hbox{Left hand side of \eqref{U_with_ell}}\\
&& \le   \big|\cQ_x(t, \eta, \gl)-\cP_x(t, \eta, \gl)\big|+ 
\big|\cP_x(t, \eta, \gl)-x_\ell(t) \wb g\big( \gs(\eta -\phi^\gp_\ell(t)), \gs^2(\gl - \phi^{\gp \gp}_\ell(t))\big) \big| \\
&& \le M(t)\Pi(t)+\sum_{k\not= \ell } A_k(t)  \gUp_{\ell, k}(t),   
\end{eqnarray*}
where the last inequality follows from \eqref{S_R_error} and \eqref{middel_step}. This shows \eqref{U_with_ell}. 
\hfill $\blacksquare$ 

\begin{theo}\label{theo:CT3S1}
Let $x(t)\in \cE_{\ep_1,\ep_3}$ for some $\ep_1, \ep_3>0$. 
Suppose that \eqref{def_sep_cond_cros} holds for some $\rho, \gt>0$, and 
\begin{equation}
\label{theo1_cond1}
2M(t)\big(\gUp(t)+\Pi(t)\big)\le \mu(t). 
\end{equation} 
Also, let $\cG_t$ and $\cG_{t, k}$ be the sets defined by \eqref{def_cGk} for some $\wt \ep_1=\wt \ep_1(t)>0$ that  satisfies    
\begin{equation}
\label{cond_ep1}
M(t)\big(\gUp(t)+\Pi(t)\big)\le 
\wt \ep_1 \le \mu(t)-M(t)\big(\gUp(t)+\Pi(t)\big).
\end{equation}
Then $\cG_t$ is the disjoint union of the non-empty sets ${\cal G}_{t, k}, 0\le k\le K$. 
\end{theo}


{\it Proof}. \; 
Let us first prove that $\mathcal{G}_t=\cup_{k=0}^K {\cal G}_{t, k}$. Since it is clear that
 $\cup_{k=0}^K {\cal G}_{t, k}\subseteq \cG_t$, it is sufficient to show that $\cG_t\subseteq \cup_{k=0}^K \cG_{t, k}$.  Consider $(\eta, \gl)\in \cG_t$ and assume $(\eta, \gl) \not \in \cG_{t, k}$ for all $k$.  Then by applying \eqref{S_R_error} and \eqref{def_upper_bounds}, we have
 \begin{eqnarray*}
\big|\cQ_x(t, \eta, \gl)\big|\hskip -0.6cm && \le \big|\cQ_x(t, \eta, \gl)-\cP_x(t, \eta, \gl)\big|+\big|\cP_x(t, \eta, \gl)\big|\\
&&\le  M(t)\Pi(t)+\Big|\sum_{k=0}^K x_k(t) \wb g\big( \gs(\eta -\phi^\gp_k(t)), \gs^2(\gl - \phi^{\gp \gp}_k(t))\big)\Big|\\
&&\le M(t)\Pi(t)+ \sum_{k=0}^K A_k(t) \gUp(t)\\
&& = M(t)\Pi(t)+ M(t)\gUp(t) \le \wt \ep_1,   
\end{eqnarray*}
so that $(\eta, \gl)\not \in \cG_t$, which is a contradiction. Thus, $(\eta, \gl)\in \cG_{t,\ell}$ for some $\ell$. Therefore, we have $\cG_t=\cup_{k=0}^K \cG_{t, k}$. 

\bigskip 

Let $Z_k, 0\le k \le K$, be the sets defined by \eqref{def_Zk}. We claim that these sets do not overlap. 
If, on the contrary, that there exists some $(t, \eta, \gl)\in Z_\ell \cap Z_{ k}$ for $\ell\not=k$, then 
\begin{eqnarray*}
&&|\phi'_k(t)-\phi'_\ell(t)|+\rho \; |\phi''_k(t)-\phi''_\ell(t)|\\
&&\le |\phi'_k(t)-\eta|+\rho |\phi''_k(t)-\gl|+|\phi'_\ell(t)-\eta|+\rho |\phi''_\ell(t)-\gl|\\
&& \le 2 \gt.    
\end{eqnarray*}
This violates the inequality in \eqref{def_sep_cond_cros}. In other words,  \eqref{def_sep_cond_cros} implies that the sets $Z_k$, for $0\le k\le K$, are disjoint. This, together with the fact that  $\cG_{t, k}=\cG_t\cap \{(\eta, \gl): \; (t, \eta, \gl)\in Z_ k\}$,  leads to 
$\cG_{t, k}, 0\le k \le K$ are disjoint.

\bigskip
Finally let us show that each $\cG_{t, \ell}$ is a non-empty set. In this regard,  we prove that $(\phi'_\ell(t), \phi''_\ell(t))\in \cG_t$. Since $\int_{-\infty}^\infty g(t)dt=1$, we have $\wb g(0, 0)=1$. By considering   
 $\eta=\phi'_\ell(t), \gl=\phi''_\ell(t)$ in \eqref{U_with_ell}, we have 
\begin{eqnarray*}
&&
\big|\cQ_x(t, \phi'_\ell(t), \phi''_\ell(t))\big|\ge \big|x_\ell(t) \wb g(0, 0)\big|- \Err_\ell(t)\\
&& =    A_\ell(t)-M(t)\Pi(t)-\sum_{k\not= \ell} A_k(t) \gUp_{\ell, k}(t)\\
&& > \mu(t)-M(t)\Pi(t)-M(t)\gUp(t)\ge \wt \ep_1, 
\end{eqnarray*}
yielding  $(\phi'_\ell(t), \phi''_\ell(t))\in \cG_t$. This, together with $(t, \phi'_\ell(t), \phi''_\ell(t))\in Z_\ell$, implies that  $(\phi'_\ell(t), \phi''_\ell(t))\in \cG_{t,\ell}$. Therefore, $\cG_{t,\ell}$ is a non-empty set. 
\hfill $\blacksquare$ 

\bigskip 
In Algorithm 1 (CT3S),  $\phi^\gp_\ell(t)$ and $x_\ell(t)$ are approximated by $\wh \eta_\ell(t)$ and $\cQ_x(t,\wh \eta_\ell, \wh \gl_\ell)$, respectively. To establish the approximation error bounds, we 
need to impose certain conditions on the window functions. Such window functions are said to be admissible in \cite{CJLL21_adpCWT_cross}, as follows.  

\begin{mdef}\label{definition2} {\rm ({\bf Admissible window functions})} \; 
A function  $g(t)$ in $L_2(\R)\cap L_1(\R)$ with $\int_\RR g(t)dt=1$ is said to be an admissible window function if it satisfies the following conditions.
\begin{itemize}
\item[{\rm (a)}] There exists a constant $C$ such that
\begin{equation}\label{inequality_g1}
|\wb g (\eta, \gl)|\le \frac  {C}{\sqrt{ |\eta|+ |\lambda|}}, \; \forall \eta, \gl\in \RR.
\end{equation}
\item[{\rm (b)}]  $|\wb g(\eta, \gl)|$ can be written as $f(|\eta|, |\gl|)$ for some function $f(\xi_1, \xi_2)$ defined on $0\le \xi_1, \xi_2< \infty$. 

\item[{\rm (c)}] There exists $b_0$ with $0<b_0<1$ and strictly decreasing non-negative continuous functions $\gb(\xi)$ and $\gga(\xi)$ on $[0, \infty)$ with $\gb(0)=1$, $\gga(0)=1$ 
such that if the function $f$, as defined in {\rm (b)}, satisfies 
\begin{equation}
\label{ineq_cond}
1-b\le f(\eta, \gl), 
\end{equation} 
for $b$ with $0\le b \le b_0$ and $\eta, \gl$, then 
 \begin{equation}
\label{cond_gb_gga}
 1-b\le \gb(\eta), \quad  1-b\le \gga(\gl).
\end{equation}  
\end{itemize}
\end{mdef}

In Section 3, we will show  that the Gaussian function in \eqref{def_g} is an admissible window function and will derive explicit expressions of the companion decreasing functions $\gb$ and $\gga$.

\begin{theo}\label{theo:CT3S2}
Let $g$ be an admissible window function for certain positive $b_0$ such that 
the condition (c) in Definition 1 holds.
Suppose that $x(t)\in \cE_{\ep_1,\ep_3}$ for some $\ep_1, \ep_3>0$ such that \eqref{def_sep_cond_cros} and  \eqref{theo1_cond1} hold and that for $0\le \ell \le K$, $2\Err_\ell(t)/A_\ell(t)\le b_0$, where $\Err_\ell(t)$ is defined by \eqref{def_Err}. Let $\cQ_x(t,\eta,\lambda)$ be the adaptive CT of $x(t)$ with $g$, and $\cG_t$ and $\cG_{t, k}$ be the sets defined by \eqref{def_cGk} for some 
$\wt \ep_1$ satisfying  \eqref{cond_ep1}. Then for $\wh \eta_\ell(t), \wh \gl_\ell(t)$ as defined in \eqref{def_max_eta},  
\begin{enumerate}
\item[{\rm (a)}] for $\ell=1, 2, \cdots, K$, 
\begin{eqnarray}
&&\label{phi_est}
|\wh\eta_{\ell}(t)-\phi_{\ell}^\gp(t)|\le \Bd_{1, \ell}:=\frac{1}{\gs(t)} \gb^{-1}\big(1-\frac {2 \; \Err_\ell(t)}{A_\ell(t)}\big), \\
&&\label{phi_est_gl}
|\wh\gl_{\ell}(t)-\phi_{\ell}^{\gp\gp}(t)|\le \Bd_{2, \ell}:=\frac{1}{\gs^2(t)} \gga^{-1}\big(1-\frac {2 \; \Err_\ell(t)}{A_\ell(t)}\big); 
\end{eqnarray}
\item[{\rm (b)}] for $\ell=0, 1, \cdots, K$, 
\begin{equation}
\label{comp_xk_est}
\big|\cQ_{x}(t, \wh \eta_\ell, \wh \gl_\ell)- x_\ell(t)\big|
\le\Bd_{3, \ell},
\end{equation}
where 
\begin{equation*}
  \Bd_{3, \ell}:=\Err_\ell(t)+2\pi I_1 A_\ell(t) \gb^{-1}\big(1-\frac {2 \; \Err_\ell(t)}{A_\ell(t)}\big)
+\pi I_2 A_\ell(t) \gga^{-1}\big(1-\frac {2 \; \Err_\ell(t)}{A_\ell(t)}\big)  
\end{equation*}
with $I_1$ and $I_2$ defined by \eqref{def_In}; 
\item[{\rm (c)}] 
 if, in addition, $g(t)\ge 0$ for $t\in \RR$, then for $\ell=0, 1, \cdots, K$,
\begin{equation}
\label{abs_IA_est}
\big| |\cQ_x(t,\wh \eta_{\ell}, \wh \gl_\ell)|-A_{\ell}(t) \big|\le \Err_\ell(t). 
\end{equation}
\end{enumerate}
\end{theo}

\bigskip 
{\bf Proof of Theorem \ref{theo:CT3S2}(a).} 
First,  letting $\eta=\phi_\ell^\gp(t), \gl=\phi_\ell^{\gp\gp}(t)$ in \eqref{U_with_ell}, we have 
\begin{equation}
\label{est_xk2}
\big|\cQ_x(t, \wh \eta_\ell, \wh \gl_\ell)\big|  \ge \big|\cQ_x(t, \phi_\ell^\gp(t), \phi_\ell^{\gp\gp}(t))\big |\ge |x_\ell(t) \wb g(0, 0) \big| - \Err_\ell(t)
=A_\ell(t)-   \Err_\ell(t), 
\end{equation}
where the first inequality above follows from the definition of $\wh \eta_\ell, \wh \gl_\ell$ and the fact $\wb g(0, 0)=1$. In addition, \eqref{U_with_ell} implies that
\begin{equation}
\label{est_xk1}
\big|\cQ_x(t, \wh \eta_\ell, \wh \gl_\ell)\big| \le 
\big|x_\ell(t) \wb g\big(\gs(\wh \eta_\ell-\phi_\ell^\gp(t)), \gs^2(\wh \gl_\ell-\phi_\ell^{\gp\gp}(t)) \big)\big|+ \Err_\ell(t). 
\end{equation}
Thus
$$
A_\ell(t) - \Err_\ell(t)\le A_\ell(t) 
\big|\wb g\big(\gs(\wh \eta_\ell-\phi_\ell^\gp(t)), \gs^2(\wh \gl_\ell-\phi_\ell^{\gp\gp}(t)) \big)\big|
+ \Err_\ell(t),  
$$
and hence, 
\begin{equation}
\label{est_xk3}
1- \frac{2\; \Err_\ell(t)}{A_\ell(t)}\le  f\big(\gs|\wh \eta_\ell-\phi_\ell^\gp(t)|, \gs^2|\wh \gl_\ell-\phi_\ell^{\gp\gp}(t)| \big).  
\end{equation}
Under the conditions \eqref{ineq_cond} and \eqref{cond_gb_gga} for $f$, we may conclude, from \eqref{est_xk3}, that
$$
1- \frac{2\; \Err_\ell(t)}{A_\ell(t)}\le \gb \big(\gs\big|\wh \eta_\ell-\phi_\ell^\gp(t)\big|\big), \quad 
1- \frac{2\; \Err_\ell(t)}{A_\ell(t)}\le \gga \big(\gs^2\big|\wh \gl_\ell-\phi_\ell^{\gp\gp}(t)\big|\big).  
$$
Hence, since $\gb(\xi)$ and $\gga(\xi)$ are decreasing functions,
$$
\gs\big| \wh \eta_\ell-\phi_\ell^\gp(t) \big|\le \gb^{-1}\big(1- \frac{2\; \Err_\ell(t)}{A_\ell(t)}\big), \;
\gs^2\big| \wh \gl_\ell-\phi_\ell^{\gp\gp}(t) \big|\le \gga^{-1}\big(1- \frac{2\; \Err_\ell(t)}{A_\ell(t)}\big),  
$$
so that \eqref{phi_est} and \eqref{phi_est_gl} hold.

\bigskip 

{\bf Proof of Theorem \ref{theo:CT3S2}(b).} 
Applying \eqref{U_with_ell},  \eqref{phi_est} and \eqref{phi_est_gl}, we have 
\begin{eqnarray*}
&&\big|\cQ_x(t, \wh \eta_\ell, \wh\gl_\ell)- x_\ell(t)\big|\le \big|\cQ_x(t, \wh \eta_\ell, \wh\gl_\ell) - x_\ell(t) \wb g\big( \gs(\wh \eta_\ell -\phi^\gp_\ell(t)), \gs^2(\wh \gl_\ell - \phi^{\gp \gp}_\ell(t))\big)\big|\\
&& \qquad +\big |x_\ell(t) \wb g\big( \gs(\wh \eta_\ell -\phi^\gp_\ell(t)), \gs^2(\wh \gl_\ell - \phi^{\gp \gp}_\ell(t))\big)
- x_\ell(t)\big|
\\
&&\le \Err_\ell(t)+A_\ell(t) \Big| \int_\RR \frac 1\gs g(\frac \tau\gs)\Big(e^{-i2\pi (\wh \eta_\ell-\phi_\ell^{\gp}(t))\tau -i\pi (\wh \gl_\ell- \phi_\ell^{\gp\gp}(t)) \tau^2}-1\Big) d\tau \Big|\\
&&\le \Err_\ell(t)+A_\ell(t) \int_\RR \Big| \frac 1\gs g(\frac \tau\gs)\Big| \; \big|2\pi (\wh \eta_\ell-\phi_\ell^{\gp}(t))\tau +\pi (\wh \gl_\ell- \phi_\ell^{\gp\gp}(t)) \tau^2 \big| d\tau \\
&&\le \Err_\ell(t)+A_\ell(t) 2\pi |\wh \eta_\ell-\phi_\ell^{\gp}(t)|\int_\RR \big| \frac 1\gs g(\frac \tau\gs) \tau \big| d\tau +\pi \big|\wh \gl_\ell- \phi_\ell^{\gp\gp}(t)\big| \int_\RR \frac 1\gs \big|g(\frac \tau\gs)\big|  \tau^2 d\tau \\
&&= \Err_\ell(t)+A_\ell(t) 2\pi I_1 \gs |\wh \eta_\ell-\phi_\ell^{\gp}(t)|+ A_\ell(t)\pi I_2 \gs^2
\big|\wh \gl_\ell- \phi_\ell^{\gp\gp}(t)\big|\\
&&\le  \Err_\ell(t)+2\pi I_1 A_\ell(t) \gb^{-1}\big(1-\frac {2 \; \tErr_\ell(t)}{A_\ell(t)}\big)+
\pi I_2 A_\ell(t) \gga^{-1}\big(1-\frac {2 \; \tErr_\ell(t)}{A_\ell(t)}\big), 
\end{eqnarray*}
so that \eqref{comp_xk_est} holds. 

\bigskip 
{\bf Proof of Theorem \ref{theo:CT3S2}(c).}  The assumptions $g(t)\ge 0$ and $\int_\RR g(t)dt=1$ lead to 
$$
|\wb g(\eta, \gl)|
\le \int_{\RR} \Big| g(\tau) e^{-i2\pi\eta\tau-i\pi \lambda \tau^2}\Big| d\tau
=\int_\RR g(\tau)d\tau
=1
$$ for any $\eta, \gl\in \R$. Thus, from \eqref{est_xk1}, we have    
\begin{equation*}
\big|\cQ_x(t, \wh \eta_\ell, \wh \gl_\ell)\big| \le 
A_\ell(t) + \Err_\ell(t). 
\end{equation*}
This, together with \eqref{est_xk2},  implies \eqref{abs_IA_est}. 
\hfill $\blacksquare$ 

\section{Adaptive chirplet transform with Gaussian window function} 

In this section we consider the adaptive chirp transform with the Gaussian function given by \eqref{def_g} 
 and provide explicit expressions for the approximation error bounds  
 $\Bd_{1, \ell}$,  $\Bd_{2, \ell}$,  $\Bd_{3, \ell}$ in Theorem 2.  

First we show that the Gaussian function $g$ defined in \eqref{def_g} is an admissible window function. 
For this $g$, 
we have (see, for example (1.41) on page 10 of  \cite{Leon_Cohen}) 
\begin{equation}
\label{g_PFT} 
\wb g(\eta, \gl)=\frac 1{\sqrt{1+i2\pi\gl}} e^{-\frac{2\pi^2 \eta ^2}{1+i2\pi \gl}}.   
\end{equation}
Thus, we may write  $|\wb g(\eta, \gl)|= f(|\eta|, |\gl|)$, with 
 \begin{equation}
 \label{def_f}
 f(\eta, \gl):=\frac 1{(1+4\pi^2\gl^2)^{1/4}} e^{-\frac{2\pi^2 \eta ^2}{1+4\pi^2 \gl^2}} \; .  
 \end{equation}
First, let us verify that   \eqref{inequality_g1} holds with $C=2^{1/4} \pi^{-1/2}$. Indeed, 
\begin{eqnarray}
\nonumber |\wb g(\eta, \gl)|\hskip -0.6cm &&\le\min\big\{\frac 1{(2\pi^2 \eta ^2)^{1/4}}, 
\frac 1{(1+4\pi^2\gl^2)^{1/4}}\big\}\\
 \nonumber &&\le \frac{\sqrt 2}{\sqrt{\big((2\pi^2 \eta ^2)^{1/4})^2+\big((1+4\pi^2\gl^2)^{1/4}\big)^2}}
\\
\label{ineq_gaussian1}&&<\frac{2^{1/4}}{\sqrt \pi \sqrt{| \eta|+\sqrt 2 |\gl|}},   
\end{eqnarray}
where the first inequality above is shown in \cite{CJLL21_adpCWT_cross}. Thus, \eqref{inequality_g1} holds with $C=2^{1/4} \pi^{-1/2}$. 
 
To show that the Gaussian function $g$ satisfies the condition (c) in Definition 1, we apply the following result in \cite{CJLL21_adpCWT_cross}. 

\begin{pro} {\rm \cite{CJLL21_adpCWT_cross}} 
Let $f(\eta, \gl)$, $\gb(\eta)$ and $\gga(\gl)$ be the functions defined by \eqref{def_f} and  
  \begin{equation}
\label{def_gga} 
\gb(\eta)=e^{-2\pi^2 \eta ^2}, \; \gga(\gl)=\frac 1{(1+4\pi^2\gl^2)^{1/4}}. 
\end{equation}
 Suppose $b$ satisfies $0\le b\le 1-e^{-1/4}$. Then 
 $1-b\le f(\eta, \gl)$ implies that
 $$
 1-b \le \gb(\eta), \; 1-b\le \gga(\gl). 
 $$
 \end{pro}
By Proposition 1, we know $g$ satisfies the condition (c) in Definition 1 with $b_0=1-e^{-1/4}$, $\gb(\eta)$ and $\gga(\gl)$ defined by \eqref{def_gga}.  Therefore, the Gaussian function $g$ is an admissible window function.  

 \bigskip 
 Next we consider the quantities $\Upsilon(t)$ and $\Upsilon_{\ell, k}(t)$ for \eqref{def_upper_bounds} when $g$ is given by \eqref{def_g}.   From \eqref{ineq_gaussian1}, we have 
\begin{equation}
\label{ineq_gaussian2}
 |\wb g(\eta, \gl)|\le \frac L{(| \eta|+\rho \; |\gl|)^{\frac 12}}. 
\end{equation}
where 
$$
L:=\frac{\max\{2^{\frac 14}, \sqrt \rho\}}{\sqrt \pi}. 
$$
If $(\eta, \gl)\not \in Z_k$, then by \eqref{ineq_gaussian2}, we have
\begin{eqnarray*}
&&\big|\wb g\big( \gs(\eta -\phi^\gp_k(t)), \gs^2(\gl - \phi^{\gp \gp}_k(t))\big)\big|\le 
\frac L{\big(\gs |\eta -\phi^\gp_k(t)|+\rho \; \gs^2|\gl - \phi^{\gp \gp}_k(t)|\big)^{\frac 12}}\\
&&\le \frac {L}{\sqrt{\gs}\min\{\sqrt \gs, 1\} \big(|\eta -\phi^\gp_k(t)|+\rho \; |\gl - \phi^{\gp \gp}_k(t)|\big)^{\frac 12}}\\
&&\le \frac {L}{\sqrt{\gs} \min\{\sqrt \gs, 1\}\sqrt \gt},  
\end{eqnarray*}
where the last inequality follows from the definition of $Z_k$ given in \eqref{def_Zk}.
 Thus we may set $\Upsilon(t)=\frac {L}{\sqrt{\gs}\min\{\sqrt \gs, 1\} \sqrt \gt}$. Since $Z_\ell$ and $Z_k$ do not overlap for $\ell\not= k$, we may simply let $\Upsilon_{\ell, k}(t)=\Upsilon(t)$. 
 For such choices of $\Upsilon(t)$ and $\Upsilon_{\ell, k}(t)$,  it is clear that \eqref{def_upper_bounds} holds. On the other hand, we claim that if $\gs^2(t)|\phi^{\gp}_\ell(t)-\phi^{\gp}_k(t)|$ is large but $\gs^4(t)|\phi^{\gp\gp}_\ell(t)-\phi^{\gp\gp}_k(t)|$ is not as large,  then $\Upsilon_{\ell, k}(t)$ could be quite small, as shown in the following. 
 
First we observe that  for  $(\eta, \gl)\in Z_\ell$, 
$$
|\eta-\phi'_k(t)|< \gt, \; |\gl-\phi''_k(t)| <  \gt/\rho. 
$$  
so that for $\ell\not=k$, we have
 \begin{eqnarray*}
 && |\eta-\phi'_\ell(t)|\ge |\phi^{\gp}_\ell(t)-\phi^{\gp}_k(t)| -|\eta-\phi'_k(t)|> |\phi^{\gp}_\ell(t)-\phi^{\gp}_k(t)|-\gt, \\
 && |\gl-\phi''_k(t)| \le |\phi^{\gp\gp}_\ell(t)-\phi^{\gp\gp}_k(t)| +|\gl-\phi''_k(t)| <|\phi^{\gp\gp}_\ell(t)-\phi^{\gp\gp}_k(t)|+\gt/\rho. 
 \end{eqnarray*}
Let $h(t)$ be the function defined by 
\begin{equation}
\label{def_h}
h(t)=\frac 1{(1+t)^{1/4}}e^{- \frac {c_0}{1+t}},  \; t\in [0, \infty), 
\end{equation}
where $c_0\ge 0$. Then if $c_0\ge \frac 14$, $h(t)$ is increasing for $1+t\le 4c_0$. Hence, if 
$$
1+4\pi^2\gs^4 \big(|\phi^{\gp\gp}_\ell(t)-\phi^{\gp\gp}_k(t)|+\gt/\rho\big)^2\le 4\cdot 2\pi^2 \gs^2 \big(|\phi^{\gp}_\ell(t)-\phi^{\gp}_k(t)|-\gt\big)^2,
$$
then for $(\eta, \gl)\in Z_\ell$, we have  
 \begin{eqnarray}
\nonumber  && \big|\wb g\big( \gs(\eta -\phi^\gp_k(t)), \gs^2(\gl - \phi^{\gp \gp}_k(t))\big)\big |\\
 \nonumber  && =\frac1{\big(1+4\pi^2\gs^4(\gl-\phi^{\gp\gp}_k(t))^2\big)^{1/4}} e^{-\frac{2\pi^2 \gs^2(\eta-\phi^{\gp}_k(t))^2}{1+4\pi^2\gs^4(\gl-\phi^{\gp\gp}_k(t))^2}}\\
 \nonumber  && \le \frac1{\big(1+4\pi^2\gs^4(\gl-\phi^{\gp\gp}_k(t))^2\big)^{1/4}} e^{-\frac{2\pi^2 \gs^2(|\phi^{\gp}_\ell(t)-\phi^{\gp}_k(t)|-\gt)^2}{1+4\pi^2\gs^4(\gl-\phi^{\gp\gp}_k(t))^2}}\\
 \label{ineq_Ups}&&\le \frac1{\big(
 1+4\pi^2\gs^4 \big(|\phi^{\gp\gp}_\ell(t)-\phi^{\gp\gp}_k(t)|+\gt/\rho\big)^2
 \big)^{1/4}} e^{-\frac{2\pi^2 \gs^2(|\phi^{\gp}_\ell(t)-\phi^{\gp}_k(t)|-\gt)^2}{1+4\pi^2\gs^4 (|\phi^{\gp\gp}_\ell(t)-\phi^{\gp\gp}_k(t)|+\gt/\rho)^2}}, 
 \end{eqnarray}
 where the last inequality follows from the increasing property of $h(t)$ for $1+t\le 4c_0$ with $c_0=2\pi^2 \gs^2 \big(|\phi^{\gp}_\ell(t)-\phi^{\gp}_k(t)|-\gt\big)^2$. Thus we may let $\Upsilon_{\ell, k}(t)$ be the quantity in \eqref{ineq_Ups}, which is very small if 
 $
 2\pi^2 \gs^2(|\phi^{\gp}_\ell(t)-\phi^{\gp}_k(t)|-\gt)^2/\big(1+4\pi^2\gs^4 (|\phi^{\gp\gp}_\ell(t)-\phi^{\gp\gp}_k(t)|+\gt/\rho)^2\big)
 $ is reasonably large. This confirms our claim.

\bigskip 
We are now ready to derive the error bounds  $\Bd_{1, \ell}$,  $\Bd_{2, \ell}$,  $\Bd_{3, \ell}$ in Theorem 2. 

\bigskip 
Observe that the inverse function $\gb^{-1}(\xi)$ of the function $\gb(\gl)$ in \eqref{def_gga} is given by 

 $$
 \gb^{-1}(\xi)=\frac1{\pi \sqrt 2} \sqrt{-\ln \xi}, \;  0<\xi <1. 
 $$
Thus, if $\frac{2 \; \Err_\ell(t)}{A_\ell(t)}\le 1-e^{-1/4}$, 
then it follows from \eqref{est_xk3} and Proposition 1 that 
$$
\sqrt{1-\frac {2 \; \Err_\ell(t)}{A_\ell(t)}}\le \gb\big(\gs(t)
|\wh\eta_{\ell}(t)-\phi_{\ell}^\gp(t)|\big); 
$$
and 
\begin{eqnarray}
\nonumber &&|\wh\eta_{\ell}(t)-\phi_{\ell}^\gp(t)|\le \Bd_{1, \ell}:=\frac{1}{\gs(t)} \gb^{-1}\big(1-\frac {2 \; \Err_\ell(t)}{A_\ell(t)}\big)\\
\label{B1_est} && = \frac{1}{\gs(t) {\pi \sqrt 2}} \sqrt{-\ln \big(1-\frac {2 \; \Err_\ell(t)}{A_\ell(t)}\big)}.
\end{eqnarray}

As to function $\gga(\gl)$ in \eqref{def_gga}, its inverse is given by
 $$
 \gga^{-1}(\xi)=\frac 1{2\pi \xi^2}\sqrt{1-\xi^4};
 $$
and hence, the error bound $\Bd_{2, \ell}$ in \eqref{phi_est_gl} can be written as 
\begin{eqnarray} 
\nonumber \Bd_{2, \ell}\hskip -0.6cm &&:=\frac{1}{\gs^2(t)} \gga^{-1}\big(1- \frac{2\; \Err_\ell(t)}{A_\ell(t)}\big)\\
\label{B2_est} &&=
 \frac 1{\gs^2(t) 2\pi \big(1-\frac{2 \; \Err_\ell(t)}{A_\ell(t)}\big)^2}\sqrt{1-\big(1-\frac{2 \; \Err_\ell(t)}{A_\ell(t)}\big)^4} \; . 
\end{eqnarray}

It also follows that the error bound $\Bd_{3, \ell}$ in \eqref{comp_xk_est} for component recovery satisfies 
\begin{equation}
\label{B3_est}
  \Bd_{3, \ell}\le \Err_\ell(t)+2 e^{1/8} I_1 \sqrt{{\Err_\ell(t)}{A_\ell(t)}}
+\frac{I_2 A_\ell(t)}{2\big(1-\frac{2 \; \Err_\ell(t)}{A_\ell(t)}\big)^2}\sqrt{1-\big(1-\frac{2 \; \Err_\ell(t)}{A_\ell(t)}\big)^4}\; . 
\end{equation}

To summarize the above derivations, we have the following theorem. 
\begin{theo}\label{theo:CT3S3}
Let $x(t)\in \cE_{\ep_1,\ep_3}$ for some $\ep_1, \ep_3>0$, such that 
\eqref{def_sep_cond_cros} and  \eqref{theo1_cond1} hold and  $2\Err_\ell(t)/A_\ell(t)\le 1-e^{-1/4}$, for $0\le \ell \le K$. 
Let $\cG_t$ and $\cG_{t, k}$ be the sets defined by \eqref{def_cGk} for some $\wt \ep_1$ that satisfies  \eqref{cond_ep1}. Let $\wh \eta_\ell(t)$ and $\wh \gl_\ell(t)$ be the functions defined by \eqref{def_max_eta}.
Then \eqref{phi_est}, \eqref{phi_est_gl} and \eqref{comp_xk_est} hold with 
$\Bd_{1, \ell}, \Bd_{2, \ell}$ and $\Bd_{3, \ell}$ given  by the quantities in \eqref{B1_est}, \eqref{B2_est} and \eqref{B3_est}, respectively. 
\end{theo}

\section{Experiments}
\begin{figure}
	\centering
	\hspace{-0.5cm}
	\begin{tabular}{cc}
		\resizebox{2.0in}{1.5in}{\includegraphics{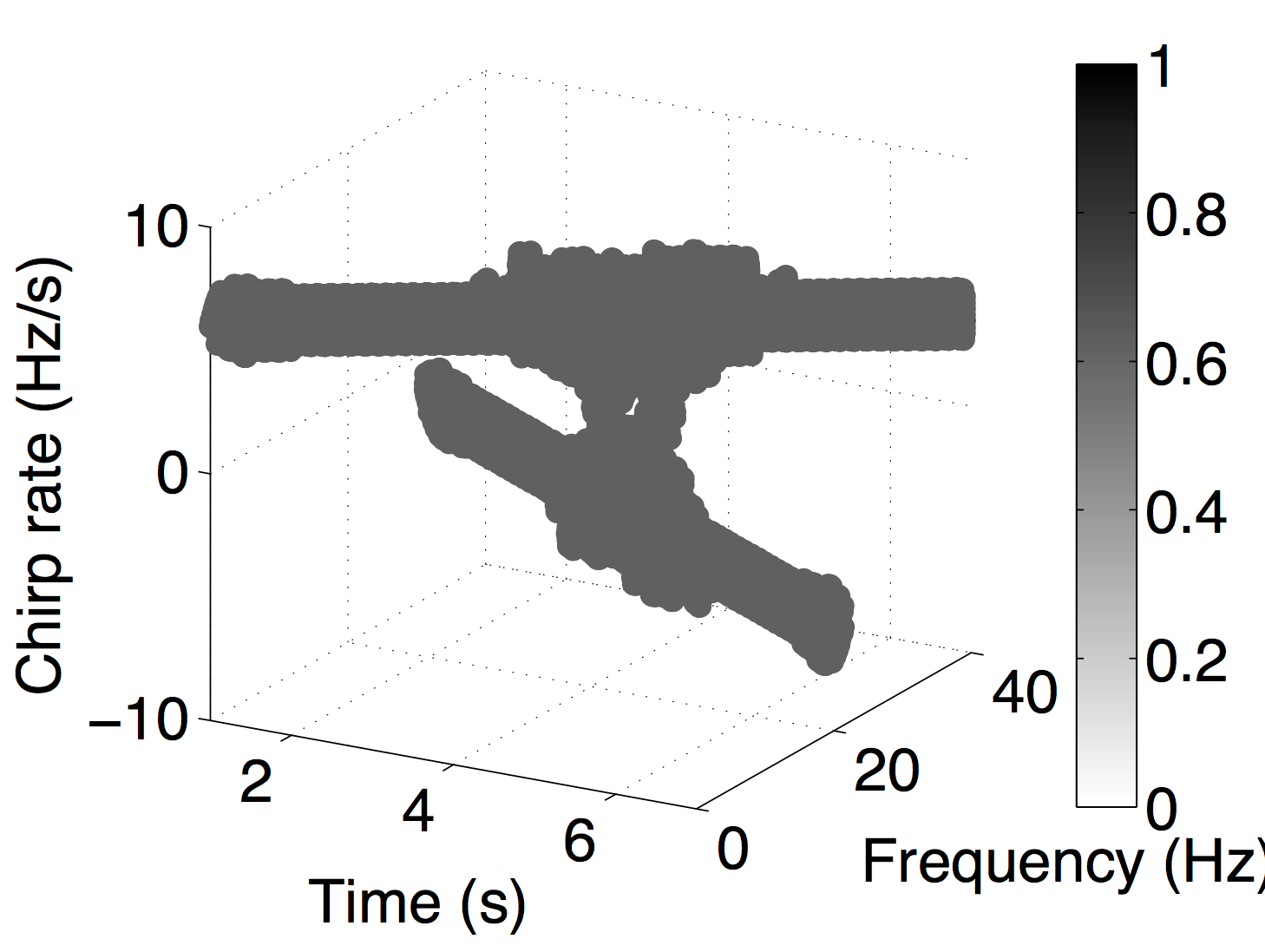}} \quad &
		\resizebox{2.0in}{1.5in}{\includegraphics{Ridge_3D_two_LFM.png}} \\
	\resizebox{2.0in}{1.5in}{\includegraphics{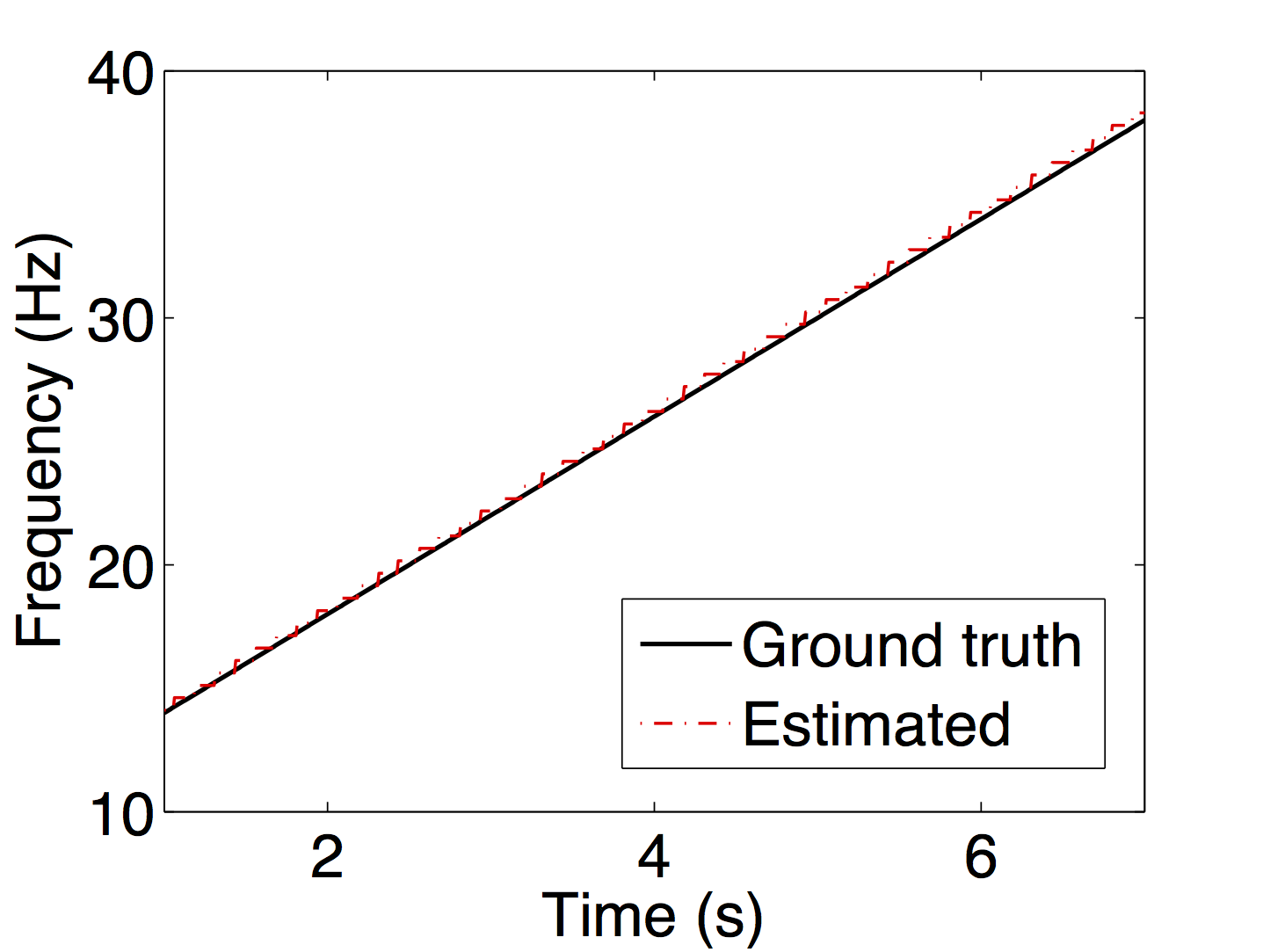}} \quad &\resizebox{2.0in}{1.5in}{\includegraphics{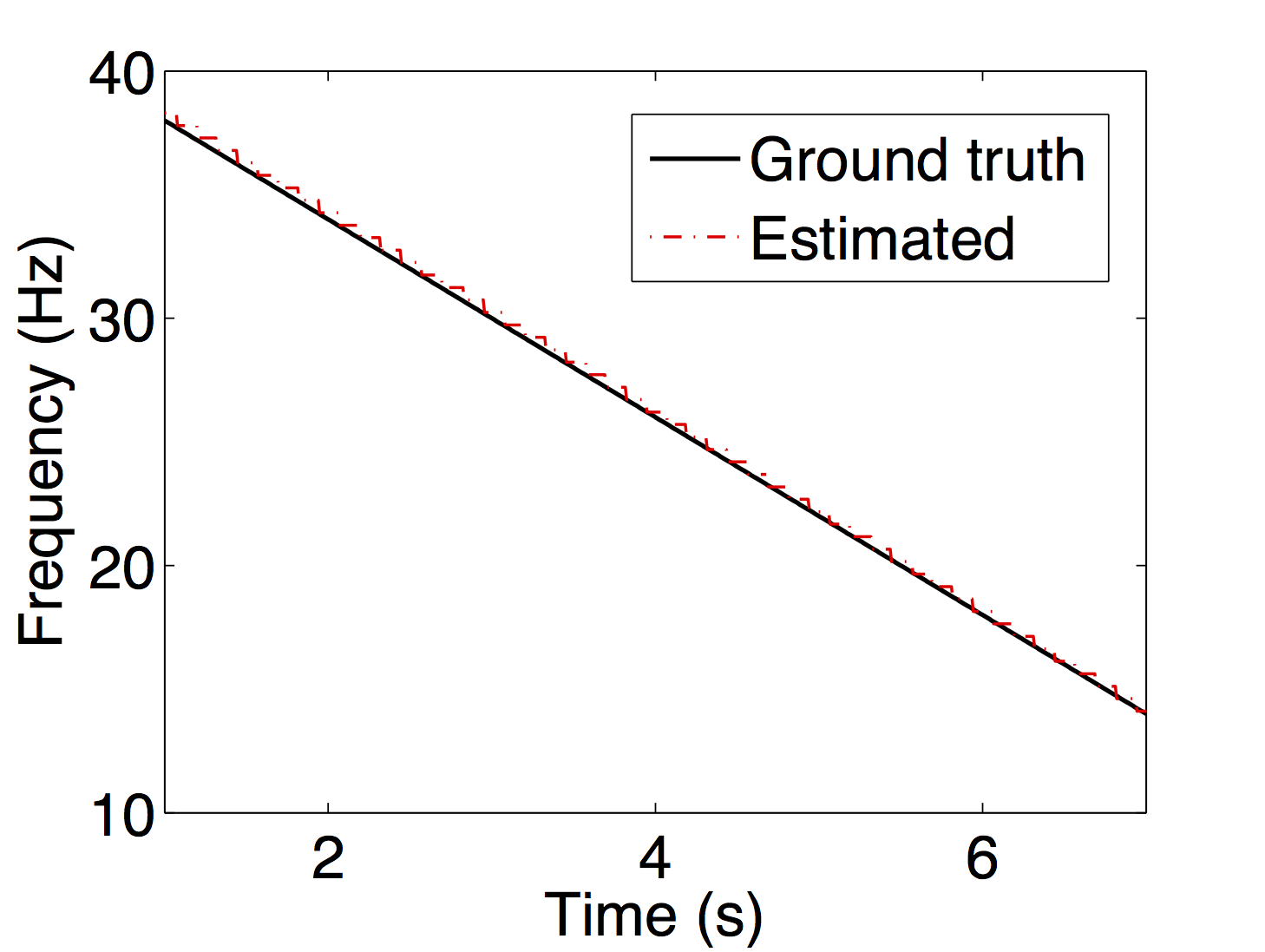}}\\
	\end{tabular}
	\caption{\small IF estimation results of the two-component LFM signal. 
		Top-right: The modulus of the CT defined by \eqref{def_ACT};
	    Top-left: The extracted ridges in the three-dimension space of CT;
        Bottom row: The ground truth and estimated IFs of Component 1 and Component 2 (from left to right).}
	\label{Fig:two-LFM}
\end{figure}

\begin{figure}
		\centering
		\hspace{-0.5cm}
		\begin{tabular}{c}
			\resizebox{4.2in}{1.5in}{\includegraphics{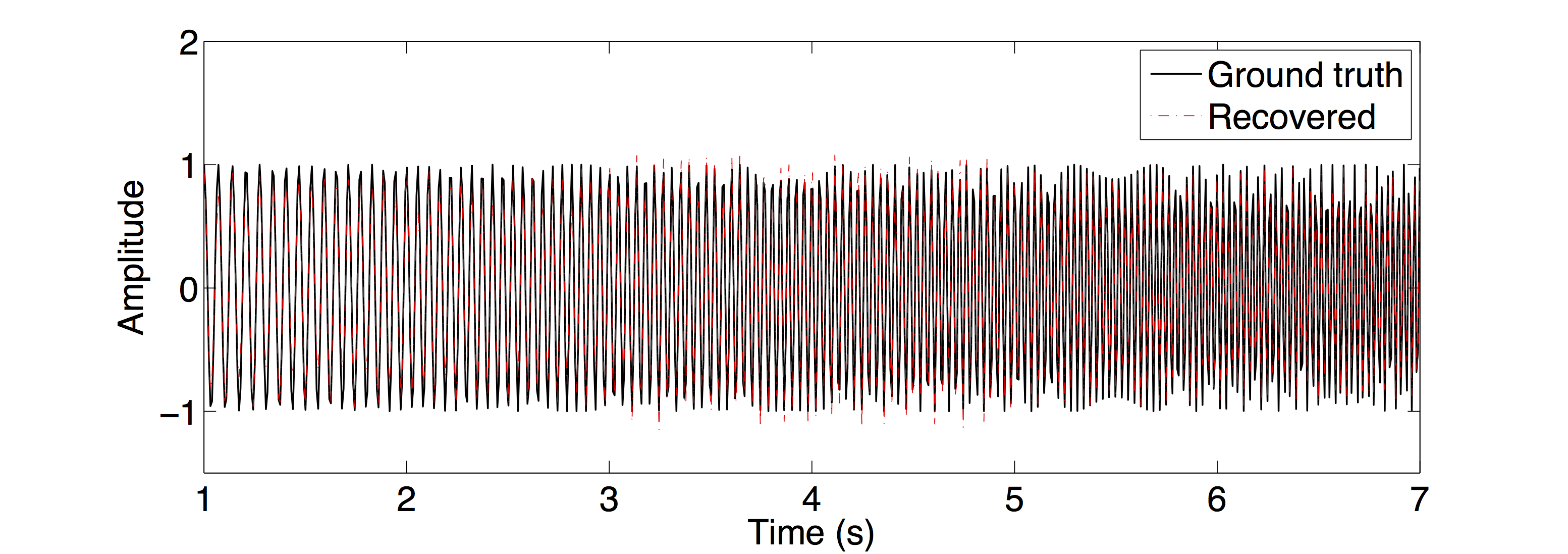}}\\ \resizebox{4.2in}{1.5in}{\includegraphics{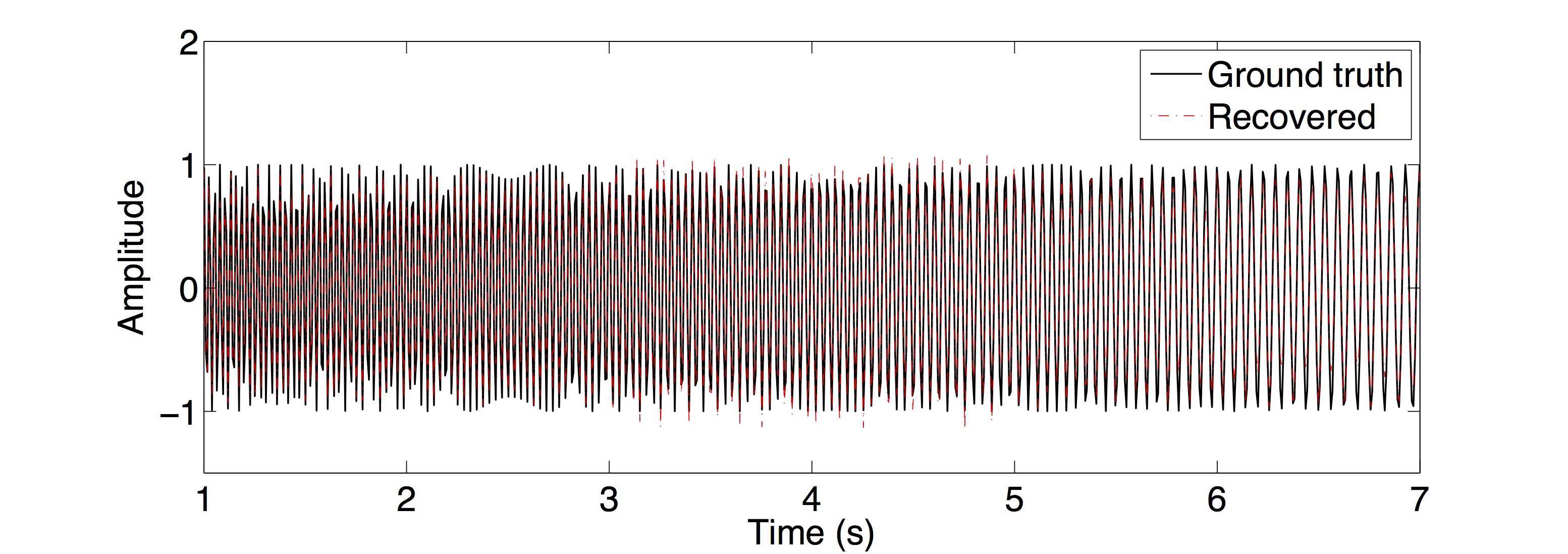}}
		\end{tabular}
		\caption{\small Recovery results of the two-component LFM signal. Top row: Component 1; Bottom row: Component 2. }
		\label{Fig:recovery of two-LFM}
	\end{figure}

In this section we present two numerical experiments to illustrate the general theorem. The interested reader is referred to \cite{LHJC22} for more experimental results for CT3S applications to IF estimation and mode retrieval.   First we consider the two-component LFM signal $x(t)$ given in 
\eqref{x_example1}, where $t$ is sampled with rate $\frac 1{128}$. In this and the other example, we let $\gs=0.15$. 
Figure \ref{Fig:two-LFM} shows the IF estimation results of the two-component LFM signal, which is based on the extracted ridges in the three-dimension (3D) space of CT defined by \eqref{def_ACT}. 
	Observe that the two-component LFM signal with crossover IFs are well separated in the three-dimensional space of time-frequency-chirp rate. 
	Furthermore, with the proposed CT3S method, we provide the sub-signal recovery results in Figure \ref{Fig:recovery of two-LFM}. Observe that the recovered modes are close to the source sub-signals.

\begin{figure}
	\centering
	\hspace{-0.5cm}
	\begin{tabular}{cc}
		\resizebox{2.0in}{1.5in}{\includegraphics{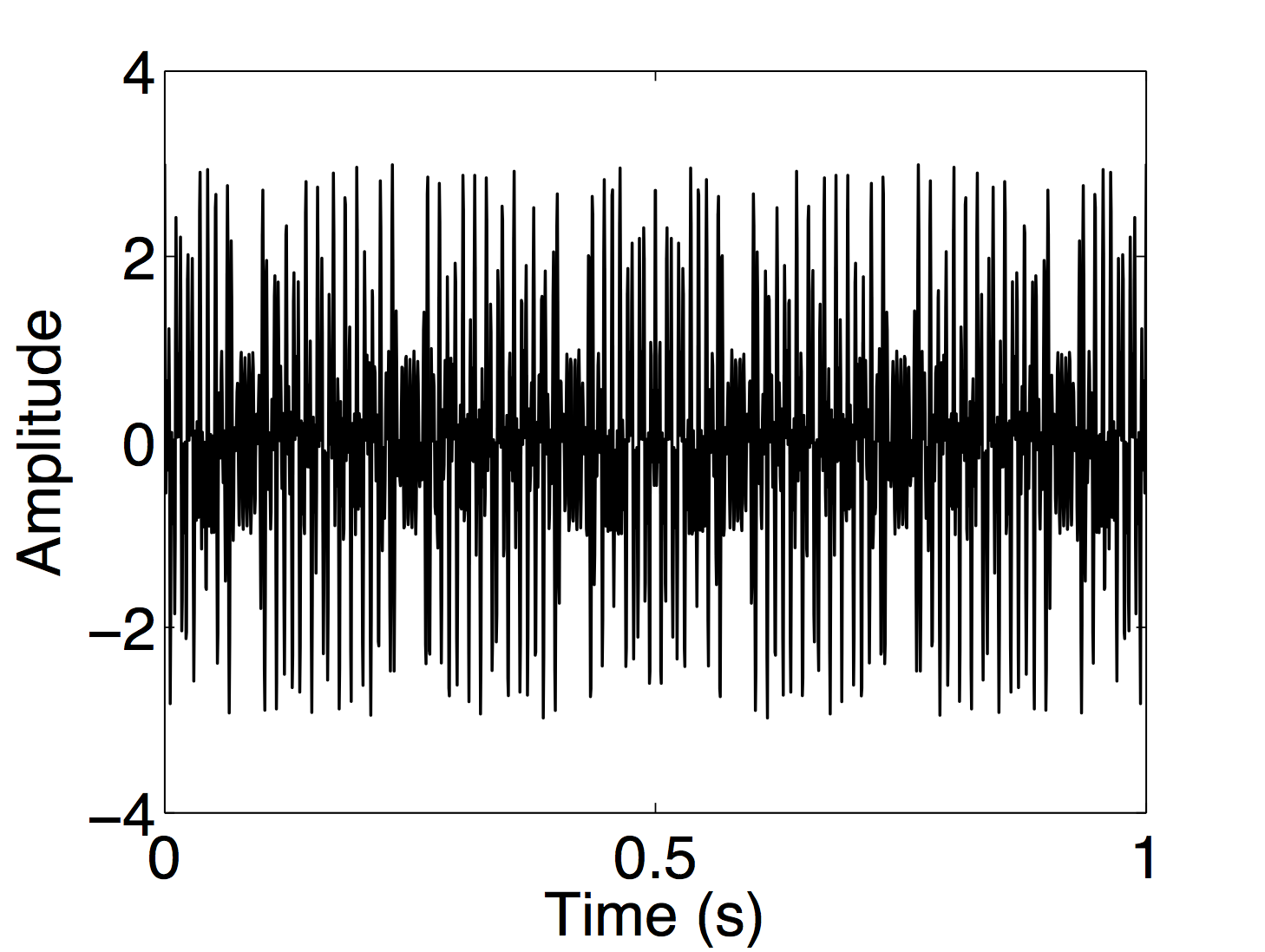}} \quad &
		\resizebox{2.0in}{1.5in}{\includegraphics{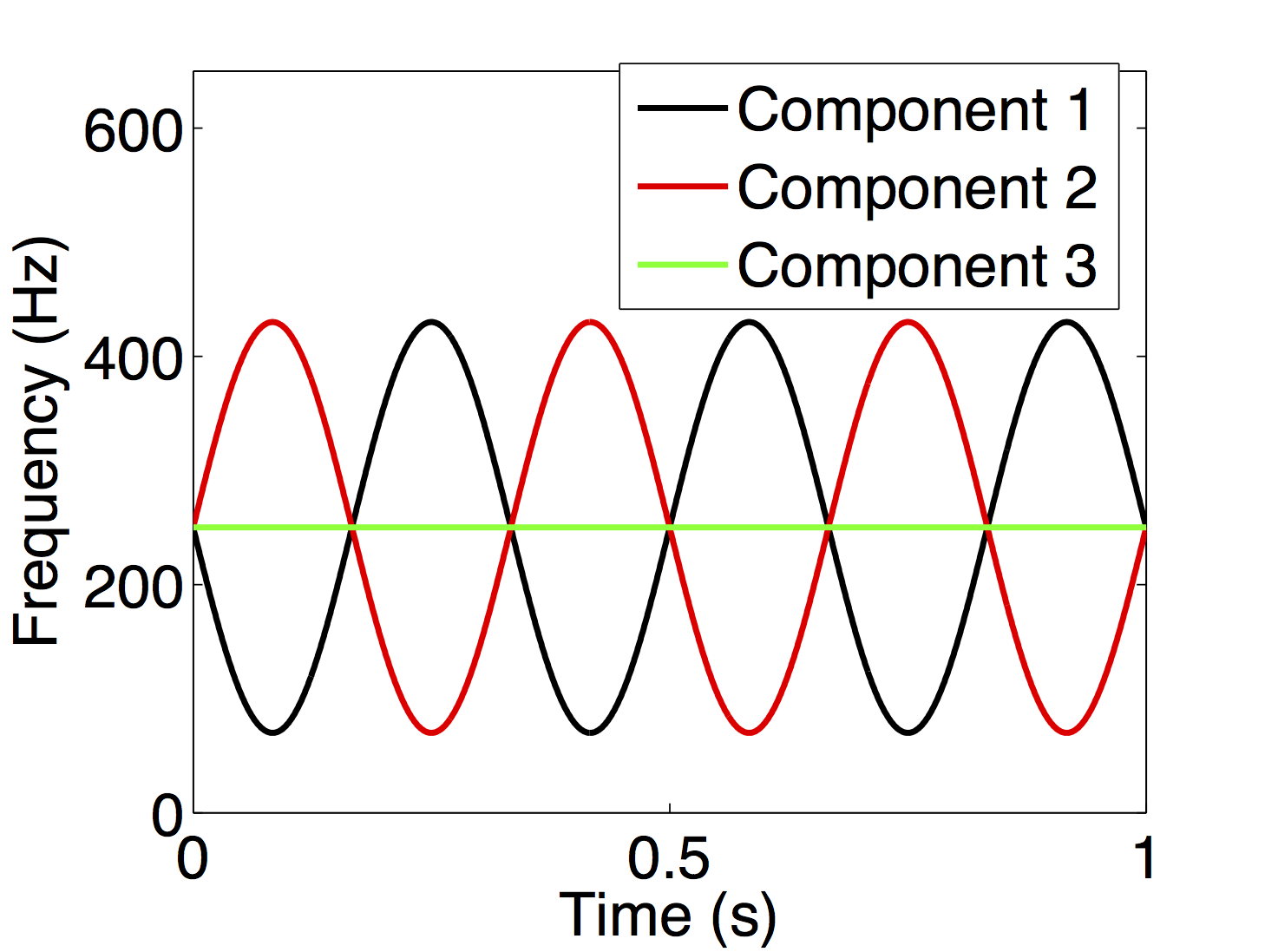}} \\
	\end{tabular}
	\caption{\small Waveform (left) and IFs (right) of the radar signal.}
	\label{radar-echoes}
\end{figure}

Next we consider the radar echoes with micro-Doppler effect overlapping in the time frequency plane \cite{radar_basic_2016}.
	The received radar echoes consist of three components, given as
	\begin{eqnarray}
	\nonumber s(t) \hskip -0.6cm &&= s_1(t) + s_2(t) + s_3(t)\\
	\nonumber && = e^{i2\pi  \left(250t - \frac{30}{\pi} \sin (6\pi t) \right)} + e^{i2\pi  \left(250t + \frac{30}{\pi} \sin (6\pi t) \right)} + e^{i500\pi t},   
	\end{eqnarray}
	where $t\in[0,1]$ and it is discretized with rate $1/2048$. Hence the IFs of $s_1(t)$, $s_2(t)$ and $s_3(t)$ are $\phi_1(t) = 250 - 180\sin(6\pi t) $, $\phi_2(t) = 250 + 180\sin(6\pi t) $ and $\phi_3(t) = 250$, respectively.  Figure \ref{radar-echoes} shows the waveform (real part of $s(t)$) and actual IFs of the radar signal.

\begin{figure}
	\centering
	\hspace{-0.5cm}
	\begin{tabular}{c}
		\resizebox{4.2in}{1.5in}{\includegraphics{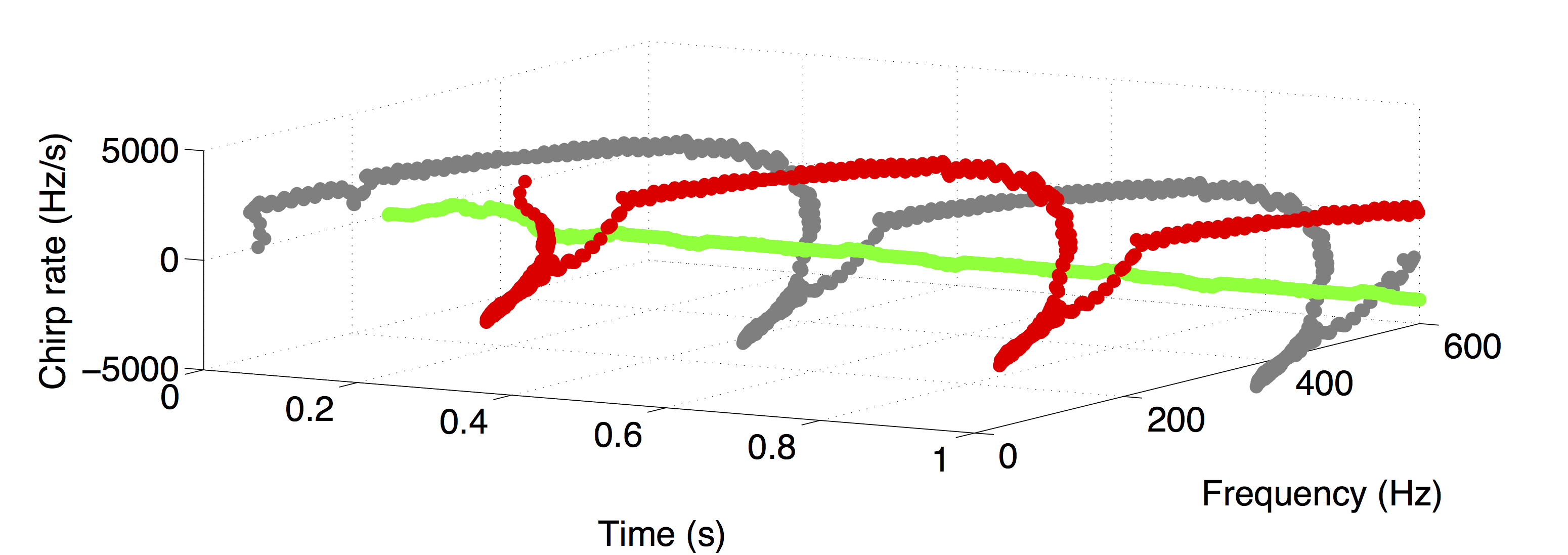}}\\
		\resizebox{4.2in}{1.5in}{\includegraphics{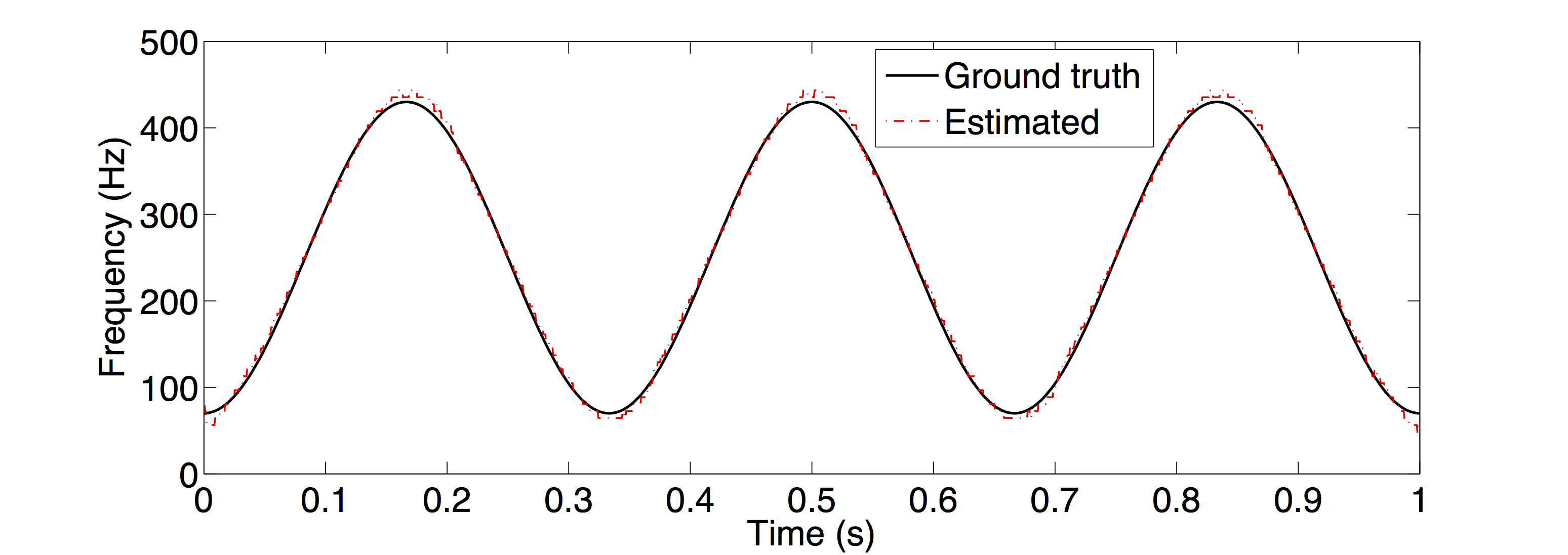}}\\
		\resizebox{4.2in}{1.5in}{\includegraphics{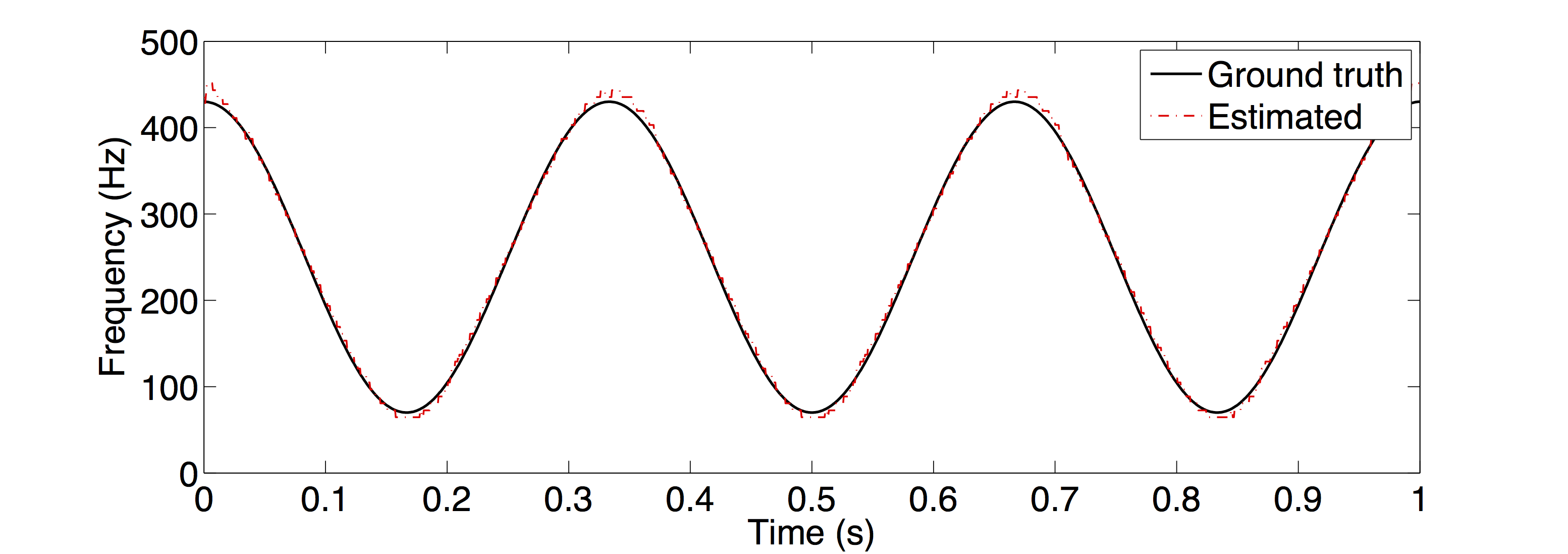}}\\
		\resizebox{4.2in}{1.5in}{\includegraphics{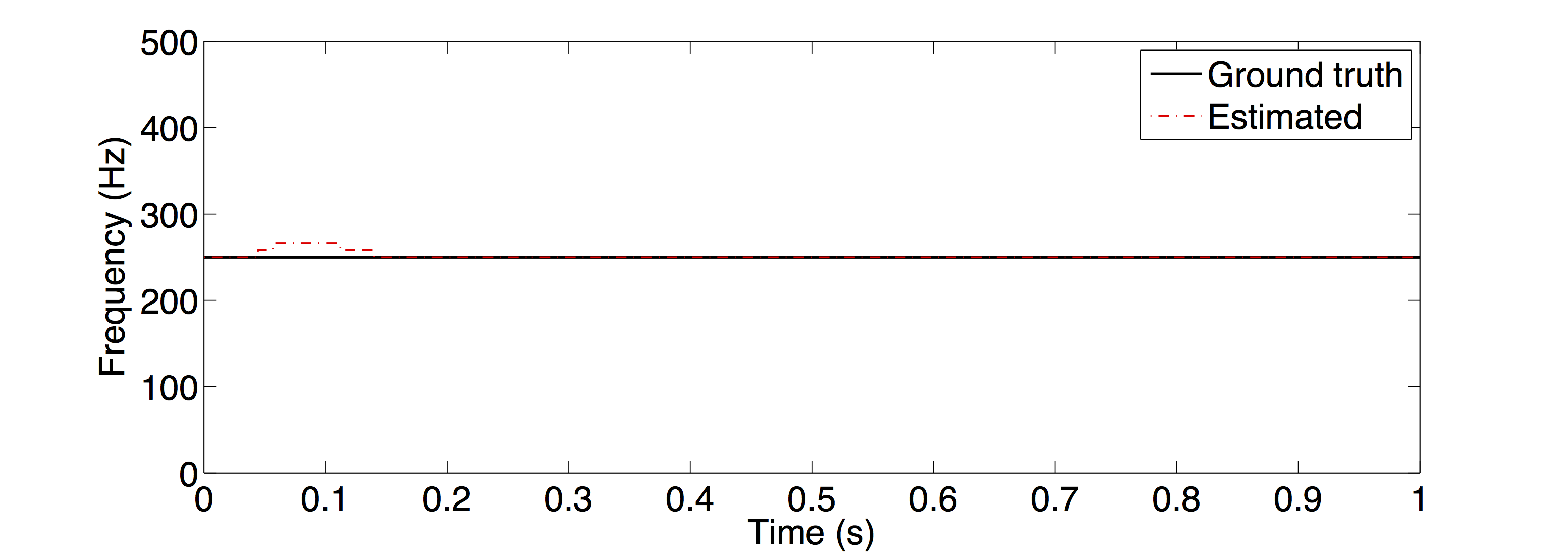}}\\
	\end{tabular}
	\caption{\small IF estimation results of the radar signal. The extracted ridges in the 3D space of ST, the estimated and ground truth IFs of $s_1(t)$, $s_2(t)$ and $s_3(t)$ (from top to bottom).}
	\label{IF-estimation-radar}
\end{figure}
	
		\begin{figure}
		\centering
		\hspace{-0.5cm}
		\begin{tabular}{c}
			\resizebox{4.2in}{1.5in}{\includegraphics{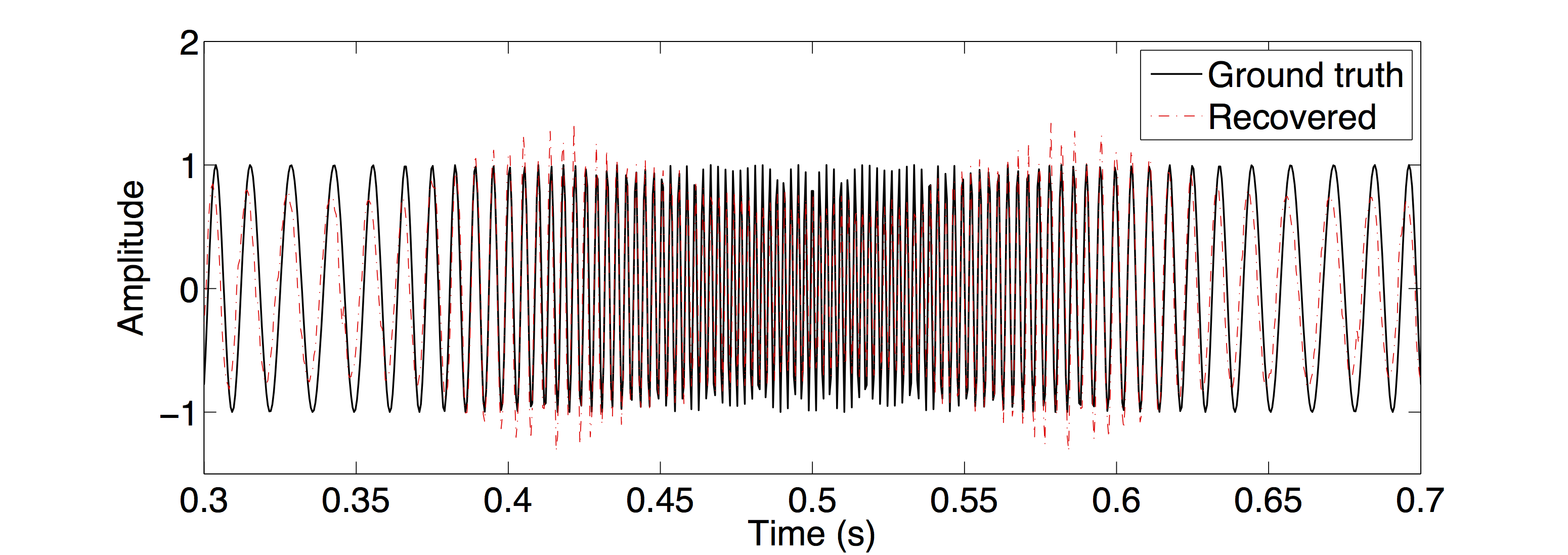}}\\
			\resizebox{4.2in}{1.5in}{\includegraphics{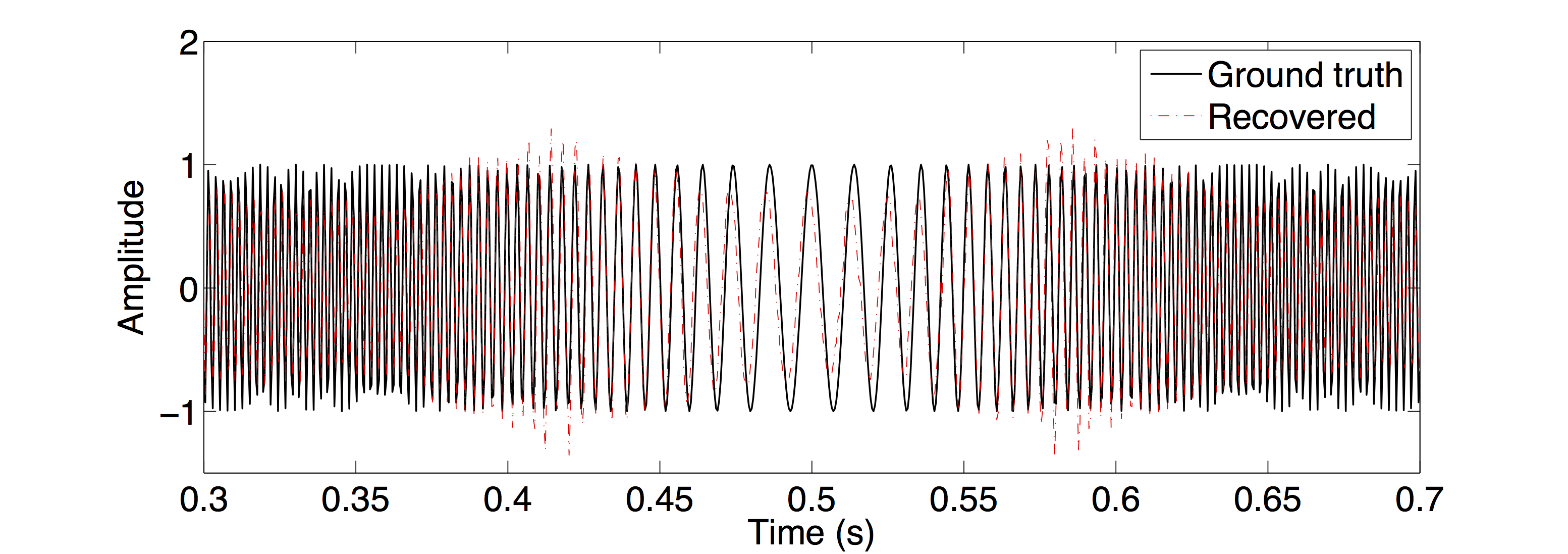}}\\
			\resizebox{4.2in}{1.5in}{\includegraphics{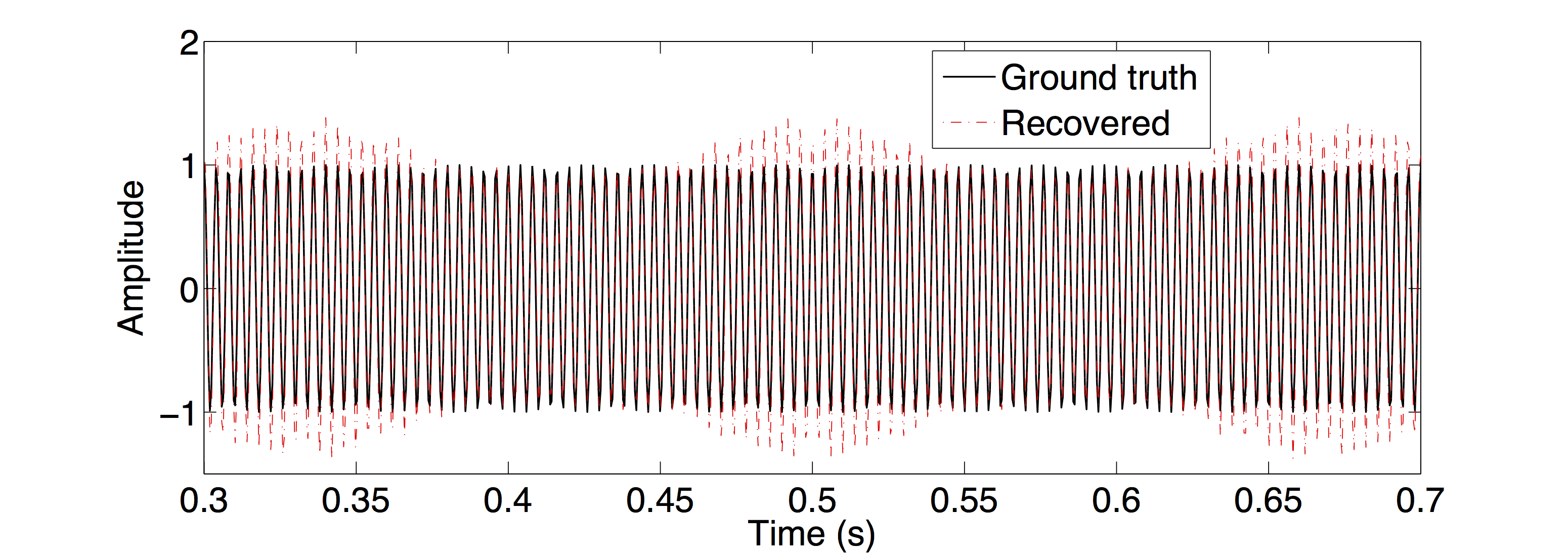}}\\
		\end{tabular}
		\caption{\small Component recovery results of the radar signal. The recovered and ground truth functions of the real parts of  $s_1(t)$, $s_2(t)$ and $s_3(t)$ (from top to bottom).}
		\label{Mode-recovery-radar}
	\end{figure}

		By using our proposed CT3S, Figure \ref{IF-estimation-radar} shows the IF estimation results. Note that the CT3S can separate the three crossover sub-signals efficiently, and estimate their IFs precisely. 
	Figure \ref{Mode-recovery-radar} shows the recovered and ground truth functions of $s_1(t)$, $s_2(t)$ and $s_3(t)$. We just plot the real parts. Meanwhile, due to the boundary effect, we just show the functions when $t\in [0.3, 0.7]$. Observe that the recovery accuracy is mainly depended on the precision of the linear chirp local approximation.

\end{document}